\documentclass[11pt]{article}
\usepackage{amssymb}
\usepackage{amsmath}
\usepackage{amsthm}
\usepackage{latexsym}
\usepackage{amsfonts}
\usepackage{graphicx}
\usepackage{graphics}

\newcommand{\dis}{\displaystyle}
\theoremstyle{plain}
\newtheorem{thm}{Theorem}[section]   
\newtheorem{prop}[thm]{Proposition}
\newtheorem{lem}[thm]{Lemma}

\theoremstyle{definition}
\newtheorem{rem}[thm]{Remark}

\newtheorem*{Proof}{Proof}

\newcommand{\bbb}[1]{\mbox{\boldmath$#1$}}

\newcommand{\el}{\ell}

\newcommand{\ra}{\;\rightarrow\;}

\newcommand{\ga}{\gamma }
\newcommand{\Ga} {{\varGamma}}

\newcommand{\de}{\delta }
\newcommand{\OO} {{\varOmega}}
\newcommand{\De} {{\varDelta}}
\newcommand{\e}{\varepsilon }
\newcommand{\f}{\varphi}

\newcommand{\zi}{\zeta }

\newcommand{\la}{\lambda }
\newcommand{\mi}{\mu }

\newcommand{\ti}{\tau }

\newcommand{\oo}{\omega}
\newcommand{\C}{\mathbb{C}}
\newcommand{\R}{\mathbb{R}}

\newcommand{\N}{\mathbb{N}}
\newcommand{\Q}{\mathbb{Q}}

\newcommand{\ssum}{\sum\limits}

\newcommand{\bz}{\bar{z}}

\newcommand{\oD}{\overline{D}}

\newcommand{\oV}{\overline{V}}
\newcommand{\oO}{\overline{\varOmega}}

\newcommand{\tr}{\tilde{r}}
\newcommand{\tp}{\tilde{p}}
\newcommand{\tq}{\tilde{q}}

\newcommand{\tf}{\widetilde{f}}
\newcommand{\tS}{\widetilde{S}}

\newcommand{\tA}{\widetilde{A}}
\newcommand{\tB}{\widetilde{B}}

\newcommand{\tz}{\tilde{z}}

\newcommand{\tit}{\tilde{t}}
\newcommand{\tiT}{\tilde{T}}

\newcommand{\cu}{{\cal{U}}}

\newcommand{\cD}{{\cal{D}}}

\newcommand{\ld}{\ldots}

\newcommand{\sm}{\smallsetminus}

\newcommand{\qb}{$\quad\blacksquare$}

 \renewcommand{\Im}{\mbox{Im}}

 \renewcommand{\Re}{\mbox{Re}}

 \newcommand{\dist}{\mbox{dist}}

\begin{document}
\title{\bf Pad\'{e} Approximants, density of rational functions in $\bbb{A^\infty(\OO)}$ and smoothness
of the integration operator}
\author{Vassili Nestoridis and Ilias Zadik}
\date{}
\maketitle
\begin{abstract}
First we establish some generic universalities for Pad\'{e}
approximants in the closure $X^\infty(\OO)$ in $A^\infty(\OO)$ of
all rational functions with poles off $\oO$, the closure taken in
$\C$ of the domain $\OO\subset\C$.\ Next we give sufficient
conditions on $\OO$ so that $X^\infty(\OO)=A^\infty(\OO)$.\ Some of
these conditions imply that, even if the boundary $\partial\OO$ of a
Jordan domain $\OO$ has infinite length, the integration operator on
$\OO$ preserves $H^\infty(\OO)$ and $A(\OO)$ as well.\ We also give
an example of a Jordan domain $\OO$ and a function $f\in A(\OO)$,
such that its antiderivative is not bounded on $\OO$.\ Finally we
restate these results for Volterra operators on the open unit disc
$D$ and we complete them by some generic results.
\end{abstract}
{\em AMS classification number}: primary 30K05, 30E10, 47G10, secondary 45P05. \vspace*{0.2cm} 
\section{Introduction}\label{sec1}
\noindent

Universality of Taylor series is a generic phenomenon, (\cite{Nes},
\cite{Ka}, \cite{Gr.Er}, \cite{B-GrEr-N-P}, \cite{M-N},
\cite{K-K-N}), where the partial sums of the Taylor expansion of a
holomorphic function $f$ approximate many functions on compact sets
outside the domain of definition $\OO$ of $f$.\ Recently the partial
sums of the Taylor expansion of $f$, which are polynomials, have
been replaced by some Pad\'{e} approximants of the Taylor series of
$f$, which are rational functions and may take the value $\infty$ as
well (\cite{Nes2}). This time the approximation is uniform on
compact sets $K$ with respect to the chordal metric $\chi$ on
$\C\cup\{\infty\}$.\ If the compact sets $K$ are disjoint from
$\oO$, then the universal function $f$ can be chosen to be smooth on
the boundary of $\OO$, that is, $f\in A^\infty(\OO)$.\ In addition
$f$ can be approximated by rational functions with poles off $\oO$
in the natural topology of $A^\infty(\OO)$; that is, $f\in
X^\infty(\OO)$, where $X^\infty(\OO)$ denotes the closure in
$A^\infty(\OO)$ of the set of rational functions with poles off
$\oO$.\ Thus, we obtain generic universalities of Pad\'{e}
approximants in $X^\infty(\OO)$; this is established in Sections
\ref{sec3} and \ref{sec4}.\ However, generic results on closed
subspaces of $A^\infty(\OO)$ may be considered as simple existence
results and therefore, they are less significant than generic
results on the whole space $A^\infty(\OO)$. That is why in Section
\ref{sec5} we give sufficient conditions such that
$X^\infty(\OO)=A^\infty(\OO)$.\ A simple form of these conditions is
that $(\oO)^0=\OO$, $\C\sm\oO$ is connected and that there exists a
constant $M<+\infty$ such that any two points of $\OO$ can be joined
in $\OO$ by a curve with length at most $M$. Under these hypotheses
polynomials are dense in $A^\infty(\OO)$. If $\OO$ is a Jordan
domain with rectifiable boundary then the above condition is
fulfilled \cite{M-N}; however, it is possible that a Jordan domain
$\OO$, whose boundary has infinite length, satisfies the previous
sufficient condition and therefore $X^\infty(\OO)=A^\infty(\OO)$.\
Such examples are all starlike Jordan domains $\OO$ where
$\partial\OO$ has infinite length.

If $\OO$ is any Jordan domain $\OO$, where $\partial\OO$ has finite
length, then it is known that the integration operator on $\OO$ is
smooth. That is, if $f\in H^\infty(\OO)$ is a bounded holomorphic
function, then its antiderivative $F$ $(F'=f$ on $\OO)$ belongs to
$A(\OO)$ and extends continuously on $\oO$. In Section \ref{sec6} we
investigate the smoothness of the integration operator $H(\OO)\ni
f\ra F(f)\in H(\OO)$ where $F'(f)=f$ on $\OO$ and $F(f)(z_0)=0$ for
some fixed point $z_0\in\OO$.\ We give sufficient conditions of the
previous type so that $F(f)\in H^\infty(\OO)$ for all $f\in
H^\infty(\OO)$, as well as $F(f)\in A(\OO)$ for all $f\in A(\OO)$.\
This may occur even if $\partial\OO$ has infinite length.\
Furthermore, we give a specific example of a Jordan domain $\OO$ and
a function $f\in A(\OO)$ so that $F(f)\notin H^\infty(\OO)$.\ This
relates to the standard singular inner function
$\exp\dfrac{z+1}{z-1}$, which has previously been used by one of the
authors (see \cite{Nes3} and [\cite{Nes4} Prop. 19]). In the above
example the boundary of $\OO$ contains only one ``bad'' point, which
can not be reached from an interior point using a curve in $\OO$
with finite length.\ Also the constructed function $f$ is almost
explicit.\ At the end of Section \ref{sec6} we reformulate the
previous results in the language of Volterra operators on the open
unit disc $D$ (see \cite{Hawi} and the references their in).

In Section \ref{sec7} we give generic versions of the results of
Section \ref{sec6}.\ For instance for any Jordan domain $\OO$ we
show that the set of functions $f\in A(\OO)$ such that $F(f)\notin
H^\infty(\OO)$ is either empty or large in the topological sense,
that is $G_\de$ and dense in $A(\OO)$ endowed with the topology of
supremum norm on $\OO$.\ We also obtain a result in this direction
for Volterra operators on the open unit disc $D$.\ Finally we show
that for all holomorphic functions $g$ in a dense subset of $H(D)$
(respectively $A(D)$), there exists $f\in A(D)$ such that
$T_g(f)\notin H^\infty(D)$, where $T_g(f)$ is the antiderivative on
$D$ of $fg'$ vanishing at 0.\ An open question is to find a complete
metric topology in the set of all Jordan domains (contained in
$\oD$), so that for the generic Jordan domain $\OO$, there exists
$f\in A(\OO)$ whose antiderivative $F$ is not bounded in $\OO$ (or
at least $F\notin A(\OO)$).
\section{Preliminaries}\label{sec2}  
\noindent

Let $\OO\subset\C$ be open. We say that a holomorphic function $f$
defined on $\OO$, belongs to $A^\infty(\OO)$ if and only if for
every $\el\in\{0,1,2,\ld\}$ the $\el$th derivative $f^{(\el)}$
extends continuously on $\oO$.

In $A^\infty(\OO)$, we consider the topology defined by the
seminorms $\dis\sup_{z\in K_n}|f^{(\el)}(z)|$, where
$\el\in\{0,1,2,\ld\}$, and $(K_n)_{n\in\N}$ is a family of compact
sets in $\oO$, such that for every compact set $L$ in $\oO$ there
exists $n\in\N$ with $L\subset K_n$. Such a family is for example
the family of the sets $\oO\cap\overline{D(0,n)}$, $n\in\N$.\ With
this topology $A^\infty(\OO)$ becomes a Fr\'{e}chet space.

Now we call $X^\infty(\OO)$, the closure in $A^\infty(\OO)$ of all
rational functions with poles off $\oO$, where the closure is taken
in $\C$.

If we consider the one point compactification
$\C\cup\{\infty\}=\widetilde{\C}$ of $\C$, then a well known metric
is the chordal metric $\chi$ on $\C\cup\{\infty\}$, where
\[
\chi(a,b)=\frac{|a-b|}{\sqrt{1+|a|^2}\sqrt{1+|b|^2}}, \ \ \text{for} \ \ a,b\in\C
\]
and $\chi(a,\infty)=\dfrac{1}{\sqrt{1+|a|^2}}$ for $a\in\C$, and
$\chi(\infty,\infty)=0$; see \cite{Ahlfors}.
\begin{prop}\label{prop2.1}
Let $K\subset\C$ be a compact set and $q=\dfrac{A}{B}$ a rational
function, where the polynomials $A,B$ do not have a common root in
$\C$.\ Then there is a sequence $q_j=\dfrac{A_j}{B_j}$, $j=1,2,\ld$
where the polynomials $A_j$ and $B_j$ have coefficients in $\Q+i\Q$
and do not have any common zero in $\C$ for all $j$, such that
$\dis\sup_{z\in K}\chi(q_j(z),q(z))\ra0$ as $j\ra+\infty$.
\end{prop}

The above proposition is well known.\ See \cite{Nes2}.

Let $\zi\in\C$ be fixed and
\[
f=\sum^\infty_{n=0}a_n(z-\zi)^n
\]
be a formal power series $(a_n=a_n(f,\zi))$.\ Often this power series is the Taylor development of a holomorphic function $f$ in a neighborhood of $\zi$.\ Let $p$ and $q$ be two non negative integers.\ The Pad\'{e} approximant $[f;p/q]_\zi(z)$ is defined to be a rational function $\phi$ regular at $\zi$ whose Taylor development with center $\zi$,
\[
\phi(z)=\sum^\infty_{n=0}b_n(z-\zi)^n,
\]
satisfies $b_n=a_n$ for all $0\le n\le p+q$ and $\phi(z)=A(z)/B(z)$, where the polynomials $A$ and $B$ satisfy
\[
\deg A\le p, \ \ \deg B\le q \ \ \text{and} \ \ B(\zi)\neq0.
\]
It is not always true that such a rational function $\phi$ exists.\ And if it exists it is not always unique.\ For $q=0$, we always have such a unique $\phi$ which is
\[
[f;p/q]_\zi(z)=\sum^p_{n=0}a_n(z-\zi)^n.
\]
For $q\ge1$ the necessary and sufficient condition for existence and
uniqueness is that the following $q\times q$ Hankel determinant is
non-zero (\cite{Th.Pade})
\[
\left|\begin{array}{lllll}
        a_{p-q+1} & a_{p-q+2} & \cdot & \cdot & a_p \\
        a_{p-q+2} & a_{p-q+3} & \cdot & \cdot & a_{p+1} \\
        \cdot & \cdot & \cdot & \cdot & \cdot \\
        \cdot & \cdot & \cdot & \cdot & \cdot \\
        a_p & a_{p+1} & \cdot & \cdot & a_{p+q-1}
      \end{array}\right|\neq0,
\]
where $a_i=0$ for $i<0$.\ If this is satisfied we write $f\in\cD_{p,q}(\zi$).\ For $f\in\cD_{p,q}(\zi)$ the Pad\'{e} approximant
\[
[f;p/q]_\zi(z)=\frac{A(f,\zi)(z)}{B(f,\zi)(z)}
\]
is given by the following Jacobi formula
\[
A(f,\zi)(z)=\left|\begin{array}{lllll}
                    (z-\zi)^qS_{p-q}(f,\zi)(z) & (z-\zi)^{q-1}S_{p-q+1}(f,\zi)(z) & \cdot & \cdot & S_p(f,\zi)(z) \\
                    a_{p-q+1} & a_{p-q+1} & \cdot & \cdot & a_{p+1} \\
                    \cdot & \cdot & \cdot & \cdot & \cdot \\
                     \cdot & \cdot & \cdot & \cdot & \cdot \\
                    a_p & a_{p+1} & \cdot & \cdot & a_{p+q}
                  \end{array}\right|,
\]
\[
B(f,\zi)(z)=\left|\begin{array}{ccccc}
                    (z-\zi)^q & (z-\zi)^{q-1} & \cdot & \cdot & 1 \\
                    a_{p-q+1} & a_{p-q+1} & \cdot & \cdot & a_{p+1} \\
                    \cdot & \cdot & \cdot & \cdot & \cdot \\
                    \cdot & \cdot & \cdot & \cdot & \cdot \\
                     a_p & a_{p+1} & \cdot & \cdot & a_{p+q}
                  \end{array}\right|,
\]
with (see \cite{Th.Pade})
\[
S_k(f,\zi)(z)=\left\{\begin{array}{ccc}
                       \ssum^k_{\nu=0}a_\nu(z-\zi)^\nu, & \text{if} & k\ge0 \\
                       0, & \text{if} & k<0.
                     \end{array}\right.
\]
If $A(f,\zi)(z)$ and $B(f,\zi)(z)$ are given by the previous Jacobi formula and they do not have a common zero in a set $K$ we write $f\in E_{p,q,\zi}(K)$.\ Equivalently
\[
|A(f,\zi)(z)|^2+|B(f,\zi)(z)|^2\neq0
\]
for all $z\in K$.\ For $K$ compact this is equivalent to the existence of a $\de>0$ such that
\[
|A(f,\zi)(z)|^2+|B(f,\zi)(z)|^2>\de
\]
for all $z\in K$.\ We will also use the following ({\cite{Th.Pade}
Th.\ 1.4.4 page 30).
\begin{prop}\label{prop2.2}
Let $\phi(z)=\dfrac{A(z)}{B(z)}$ be a rational function, where the
polynomials $A$ and $B$ do not have any common zero in $\C$.\ Let
$\deg A(z)=k$ and $\deg B(z)=\la$.\ Then for every $\zi\in\C$ such
that $B(\zi)\neq0$ we have the following:
\[
\phi\in\cD_{k,\la}(\zi),
\]
\[
\phi\in\cD_{p,\la}(\zi) \ \ \text{for all} \ \ p>k,
\]
\[
\phi\in\cD_{k,q}(\zi) \ \ \text{for all} \ \ q>\la.
\]
In all these cases $\phi$ coincides with its corresponding Pad\'{e}
approximant, that is,\linebreak $[\phi;k/\la]_\zi(z)\equiv\phi(z)$
and $[\phi;p/\la]_\zi(z)\equiv\phi(z)$ for all $p>k$ and
$[\phi;k/q]_\zi(z)\equiv\phi(z)$ for $q>\la$.
\end{prop}

The Mo\"{e}bius function $z\ra\dfrac{z+1}{z-1}$, maps every orthogonal circle to the real axis that passes through 1, to a line parallel to the imaginary axis.

Thus, as $z$ varies in such a circle, $\Re\Big(\dfrac{z+1}{z-1}\Big)$ remains constant.\ This yields that $\Big|\exp\Big(\dfrac{z+1}{z-1}\Big)\Big|$ remains constant too.

More specifically, it can be checked that the unit circle is mapped into the unit circle through the mapping $z\ra\exp\Big(\dfrac{z+1}{z-1}\Big)$.

Now consider the mapping $g$ defined on the set
$\{z\in\C\,|\;\Re(z)\le1\}$
\[
g(z)=\left\{\begin{array}{cc}
              (z-1)\exp\Big(\dfrac{z+1}{z-1}\Big), & z\neq1 \\
              0, & z=1.
            \end{array}\right.
\]
Then the following proposition holds;
\begin{prop}\label{prop2.3}
There exist a Jordan domain $V$, subset of the set $S=\{z\in\C\mid
\Re(z)\le1,\Im(z)\ge0\}$, with the following
properties\vspace*{-0.2cm}
\begin{enumerate}
\item[(i)] $V$ is contained in a set bounded from two arcs that belong in $S$ and are arcs of circles orthogonal to the real axis, passing through 1. \vspace*{-0.2cm}
\item[(ii)] $V$ contains an open arc of the unit circle that ends at 1 and $1\in\partial V$.\vspace*{-0.2cm}
\item[(iii)] The function $g$ defined above is one-to-one in $V$.\vspace*{-0.2cm}
\item[(iv)] The function $\dfrac{1}{\Big\{\exp\dfrac{z+1}{z-1}\Big\}\log(1-z)}$ belong to $A(V)$, which means
that it is continuous on $\oV$ and holomorphic on $V$.
\end{enumerate}
\end{prop}
\begin{Proof}
Consider the arc of the unit circle
$A_n=\Big\{e^{it}:\dfrac{1}{n}\le t\le\dfrac{\pi}{2}\Big\}$ for
$n=1,2,\ld\;.$

Because
\[
g'(z)=\exp\bigg(\frac{z+1}{z-1}\bigg)\frac{z-3}{z-1}\neq0,  \ \ \text{for} \ \ z\in S\sm\{1\},
\]
we get that for every $z\in A_n$, there exists $r=r(z)>0$ such that
$g\mid_{D(z,r)}$ is one-to-one, where $D(z,r)\subset S$.

Thus, because $A_n$ is compact, there are $z_1,z_2,\ld,z_{m_n}\in
A_n$ and $r_1,\ld,r_{m_n}>0$, $m_n\in\N$ such that
$A_n\subset\bigcup\limits^{m_n}_{i=1}D(z_i,r_i)\subset S$, and for
every $i=1,2,\ld,m_n$, $g\mid_{D(z_i,r_i)}$ is one-to-one and
$\dfrac{\pi}{2}\ge\arg(z_1)>\arg(z_2)>\cdots>\arg(z_{m_n})\ge\dfrac{1}{n}$.

Take now $w_i\in D(z_i,r_i)\cap D(z_{i+1},r_{i+1})$, for
$i=1,2,\ld,m_n-1$.\ Define $V_{i,\e}=\{e^{it}(w-1)+1\mid\,|w|=1$,
$\arg w_{i+1}\le\arg w\le\arg w_i$ and $|t|<\e\}$,
$i=1,2,\ld,m_n-1$, where $\e>0$ is small enough such that
$V_{i,\e}\subsetneq D(z_{i+1},r_{i+1})$. Fix $i\in\{1,\ld,m_n-1\}$
and denote $V_i$ the set $V_{i,\e}$ for the previous $\e$ depending
on $i$.

Now we claim that there are two arcs of circles orthogonal to the
$x$-axis, that pass through 1, one with radius equal to $a_i<1$, and
one with radius $b_i>1$ with centers $0<1-a_i<1$ and $1-b_i<0$, such
that, if we call $W_{a_ib_i}$ the set of the points between these
arcs with argument in $(0,\arg w_1)$ $g$ satisfies the following:

If $z\in W_{a_i,b_i}\cap V_{i}$ and $\tz\in W_{a_ib_i}$ with
$g(z)=g(\tz)$, then $z=\tz$.

For $a_i<1$ and $b_i>1$, which will be determined later on, we
consider $z\in W_{a_i,b_i}\cap V_i$ and $\tz\in W_{a_i,b_i}$.

If $g(z)=g(\bz)$, then
$\dfrac{|\tz-1|}{|z-1|}=e^{\text{Re}\big(\frac{z+1}{z-1}\big)-\text{Re}\big(\frac{\tz+1}{\tz-1}\big)}$.

Now, because $z,\tz\in W_{a_i,b_i}$ and the fact that, if a complex
number $t$ belongs to a circle orthogonal to the real axis and
passes through 1 of radius $r>0$, then
$\Re\Big(\dfrac{t+1}{t-1}\Big)=2+\dfrac{1}{r}$; It follows that
$e^{\frac{1}{b_i}-\frac{1}{a_i}}\le\dfrac{|\tz-1|}{|z-1|}\le
e^{\frac{1}{a_i}-\frac{1}{b_i}}$.

Now if we choose $a_i$ and $b_i$ close enough, this will make $\tz$
to be inside $D(z_{i+1}r_{i+1})$, yielding $g(z)=g(\tz)$, where
$z,\tz\in D(z_{i+1},r_{i+1})$ and thus $z=\tz$, because $g$ is one
to one in $D(z_{i+1},r_{i+1})$.

By choosing the pairs $(a_i,b_i)$, $i=1,2,\ld,m_n-2$ to satisfy also
$W_{a_{i+1},b_{i+1}}\subset W_{a_i,b_i}$ we get that $g$ is
one-to-one on $\bigcup\limits^{m_n-2}_{i=1}(W_{a_i,b_i}\cap
V_i)=S_n$. Moreover, if $z\in S_n$ and $w\in
W_{a_{m_n-2}b_{m_n-2}}\cup S_n$ satisfy $g(z)=g(w)$, then $z=w$.

Carrying this procedure as $n$ goes to infinity by taking the union
of $S_n$, we set $V=\bigcup^\infty_{n=1}S_n$ and one can verify that
$V$ satisfies all requirements. Especially the standard singular
inner function $\exp\dfrac{z+1}{z-1}$ is far from $\infty$ and 0 on
$V$.\\
Thus, we have $\dfrac{1}{\Big[\exp\dfrac{z+1}{z-1}\Big]\log(1-z)}\in
A(V)$. \qb
\end{Proof}
\begin{lem}\label{lem2.4}
Let $h:[0,t_0]\ra\C$ be a continuous function on $(0,t_0]$, where
$t_0>0$. We assume that $\dis\lim_{t\ra0}\arg(h(t))=c\in\R$ and
$\int\limits^{t_0}_{0^+}|h(t)|dt=+\infty$. Then
$\Big|\int\limits^{t_0}_{0^+}h(t)dt\Big|$ equals also $+\infty$.
\end{lem}
\begin{Proof}
Let $t_1>0$ be such that $t_1<t_0$ and $|\arg(h(t))-s|<\dfrac{\pi}{3}$ for all $0<t<t_1$.

Then for $\tit$, $0<\tit<t_1$ it holds
\begin{align*}
\bigg|\int^{t_1}_{\tit}h(t)dt\bigg|=&\,\bigg|\int^{t_1}_{\tit}h(t)e^{is}dt\bigg|=\bigg|
\int^{t_1}_{\tit}\mid h(t)\mid e^{i(\arg(h(t))-s)}dt\bigg|\\
&\ge\Re\bigg(\int^{t_1}_{\tit}\mid h(t)\mid e^{i(\arg(h(t))-s)}dt\bigg)\\
=\,&\int^{t_1}_{\tit}\mid h(t)\mid \cos(\arg(h(t)-s))dt
\end{align*}
which is bigger than $\dfrac{1}{2}\int\limits^{t_1}_{\tit}\mid h(t)\mid dt$.\ Since $\int\limits^{t_0}_{0^+}\mid h(t)\mid dt=+\infty$, it follows easily that $\Big|\int\limits^{t_1}_{0^+}h(t)dt\Big|=+\infty$.

Moreover, the last implies that, $\Big|\int\limits^{t_0}_{0^+}h(t)dt\Big|=+\infty$, because $h$ is continuous on $(0,t_0]$. The proof is complete. \qb
\end{Proof}
\section{Smooth Universal Pad\'{e} Approximants}\label{sec3}
\noindent

For the definitions of $X^\infty(\OO)$ and the notion of the
Pad\'{e} Approximants we refer to \S\:2 and we state the following.
\begin{thm}\label{thm3.1}
Let $F\subset\N\times\N$ be a set that contains a sequence $(\tp_n,\tq_n)$, $n=1,2,\ld$, such that $\tp_n\ra+\infty$ and $\tq_n\ra+\infty$ and let $\OO\subseteq\C$ an open set.\ Let $L,\De\subset\C$ be compact sets inside $\oO$ and $K$ a compact set in $\C$ such that $K\cap\oO=\emptyset$.

Then there exists $f\in X^\infty(\OO)$ such that: for every rational
function $h$ there exists a sequence $(p_n,q_n)\in F$ $(n=1,2,\ld)$
with the following properties: \vspace*{-0.2cm}
\begin{enumerate}
\item[(i)] $f\in D_{p_n,q_n}(\zi)\cap E_{p_n,q_n,\zi}(K\cup\De)$, for every $\zi\in L$. \vspace*{-0.2cm}
\item[(ii)] For every $\el\in\N$, $\dis\sup_{\zi\in L}\dis\sup_{z\in\De}\big|\big[f;p_n/q_n\big]^{(\el)}_\zi(z)-f^{(\el)}(z)\Big|\ra0$, as $n\ra+\infty$.\vspace*{-0.2cm}
\item[(iii)] $\dis\sup_{\zi\in L}\,\dis\sup_{z\in K}\chi\big(\big[f;p_n/q_n\big]_\zi(z),h(z)\big)\ra0$, as $n\ra+\infty$.
\end{enumerate}
The set of such functions $f\in X^\infty(\OO)$ is dense and $G_\de$
in $X^\infty(\OO)$.
\end{thm}
\begin{Proof}
Let $(f_i)_{j\in\N}
$ be an enumeration of the rational functions with coefficients of the numerator and the denominator from $\Q+i\Q$.

We name $\cu$ the set of all functions in $X^\infty(\OO)$ that
satisfy the properties (i), (ii) and (iii), and we will prove that
$\cu$ is a $G_\de$-dense in the $X^\infty(\OO)$-topology and
therefore, $\cu\neq\emptyset$.

For $j,s\in\N^\ast$ and $(p,q)\in F$ we define:
\begin{align*}
E(j,p,q,s)=&\,\Big\{f\in X^\infty(\OO)\mid f\in D_{p,q}(\zi)\cap E_{p,q,\zi}(K)\ \ \text{for all} \ \ \zi\in L \ \ \text{and} \\
&\dis\sup_{\zi\in L}\,\dis\sup_{z\in K}\chi([f;p/q]_\zi(z),f_j(z))<\dfrac{1}{s}\Big\}
\end{align*}
and,
\begin{align*}
&T(p,q,s)=\{f\in X^\infty(\OO)\mid f\in D_{p,q}(\zi)\cap E_{p,q,\zi}(\De)\;
\text{for all}\; \zi\in L \\
&\text{and} \; \sup_{\zi\in L}\,\sup_{z\in\De}\Big|\Big[f;p/q\Big]^{(\el)}_\zi(z)-f^{(\el)}(z)\Big|<\dfrac{1}{s}
\;\;\text{for}\; \el=0,1,\ld,s\Big\}.
\end{align*}
Proposition \ref{prop2.1} and the definition of $X^\infty(\OO)$ imply that
\[
U=\bigcap^\infty_{j,s=1}\,\bigcup_{(p,q)\in F}(E(j,p,q,s)\cap T(p.q.s)).
\]
To prove that $U$ is a $G_\de$-dense in the
$X^\infty(\OO)$-topology, it is enough to prove that for every
$j,s=1,2,\ld$ and $(p,q)\in F$ the sets $E(j,p,q,s)$, $T(p,q,s)$ are
open in $X^\infty(\OO)$ and that for every $j$ and $s$ inside
$\N^\ast$, the set $\bigcup\limits_{(p,q)\in F}(E(j,p,q,s)\cap
T(p,q,s))$ is dense in $X^\infty(\OO)$.

Now let $j,s\in\N^\ast$ and $(p,q)\in F$. We first prove that the
set $E(j,p,q,s)$ is open in $X^\infty(\OO)$.\ Indeed, let $f\in
E(j,p,q,s)$ and let $g\in X^\infty(\OO)$ be such that,
\begin{eqnarray}
\sup_{z\in L}\mid f^{(m)}(z)-g^{(m)}(z)\mid<a \ \  \text{for} \ \  m=0,1,2,\ld,p+q+1  \label{eq1}
\end{eqnarray}

The number $a>0$ will be determined later on. It is enough to prove
that if $a$ is small enough then $g\in E(j,p,q,s)$.

The Hankel determinants defining $D_{p,q}(\zi)$ for $f$ depend
continuously on $\zi\in L$; thus, there exists $\de>0$ such that the
absolute values of the corresponding Hankel determinants are greater
than $\de>0$, for every $\zi\in L$, because $f\in D_{p,q}(\zi)$ for
every $\zi\in L$ and because $L$ is compact.

From (\ref{eq1}) we can control the first $p+q+1$ Taylor
coefficients of $g$ and by making $a>0$ small enough one can get the
Hankel determinants that define $D_{p,q}(\zi)$ to have absolute
value at least $\de/2>0$.

Therefore, $g$ will belong in $D_{p,q}(\zi)$ for every $\zi\in L$.
Now we consider the Pad\'{e} approximants of $f,g$ according to the
Jacobi formula (see preliminaries)
\[
[f;p/q]_\zi(z)=\frac{A(f,\zi)(z)}{B(f,\zi)(z)} \ \ \text{and} \ \ [g;p/q]_\zi(z)=\frac{A(g,\zi)(z)}{B(g,\zi)(z)}.
\]
Now $\mid A(f,\zi)(z)\mid^2+\mid B(f,\zi)(z)\mid^2$ vary
continuously with respect to $(z,\zi)\in K\times L$, because of the
Jacobi formula. So, there is a $\de'>0$ such that:
\[
\mid A(f,\zi)(z)\mid^2+\mid B(f,\zi)(z)\mid^2\ge\de', \ \ \text{for all} \ \ \zi\in L \ \ \text{and} \ \ z\in K.
\]
Now again from the Jacobi formula, if $a$ is small enough, one gets:
\[
\mid A(g,\zi)(z)\mid^2+\mid B(g,\zi)(z)\mid^2\ge\de'/2, \ \ \text{for all} \ \ \zi\in L \ \ \text{and} \ \ z\in K.
\]
This yields that $g\in E_{p,q,\zi}(K)$ for every $\zi\in L$. For the
rest it is enough to show that if $a$ is small enough then
$\dis\sup_{\zi\in L}\,\dis\sup_{z\in
K}\chi([g;p/q]_\zi(z),[f;p/q]_\zi(z))$ can become less than
$\dfrac{1}{s}-\dis\sup_{\zi\in L}\,\dis\sup_{z\in
K}\chi([f;p/q]_\zi(z),f_j(z))\equiv\ga>0$. By taking $a$ small as
before we have that $\mid A(f,\zi)(z)\mid^2+\mid
B(f,\zi)(z)\mid^2>\de'$ and for every $\zi\in L$ and $z\in K$ we
have $\mid A(g,\zi)(z)\mid^2+\mid B(g,\zi)(z)\mid^2>\de'/2$.

It follows that
\begin{align*}
\chi&([f;p/q]_\zi(x),[g;p/q]_\zi(z))\\
=&\,\frac{\mid A(f,\zi)(z)B(g,\zi)(z)-A(g,\zi)(z)B(f,\zi)(z)\mid}
{\sqrt{\mid A(f,\zi)(z)\mid^2+\mid B(f,\zi)(z)\mid^2}
\sqrt{\mid A(g,\zi)(z)\mid^2+\mid B(g,\zi)(z)\mid^2}} \\
&\le\frac{\sqrt{2}}{\de'}\mid A(f,\zi)(z)B(g,\zi)(z)-A(g,\zi)(z)B(f,\zi)(z)\mid
\end{align*}
for all $\zi\in L$ and $z\in K$, which easily yields the result, because the last expression can become as small as we want to, uniformly for all $\zi\in L$, $z\in K$.\ Thus, we proved that $E(j,p,q,s)$ is open.

Next, we prove that $T(p,q,s)$ is also open in $X^\infty(\OO)$. Let $f$ be a function inside $T(p,q,s)$, $L'$ be a compact set inside $\oO$ such that $L'\supset L\cup\De$ and let $g$ be a function inside $X^\infty(\OO)$ such that: $\dis\sup_{z\in L'}\mid f^{(m)}(z)-g^{(m)}(z)\mid<a$, for $m=0,1,2,\ld,\max(s,p+q+1)$, where $a>0$ will be determined later on.

In the same way as before one deduce that by making ``$a$'' small
enough it follows $g\in D_{p,q}(\zi)\cap E_{p,q,\zi}(\De)$, for all
$\zi\in L$.

Now $f(z)\in\C$, for each $z\in\De$.\ It follows that for all $\zi\in L$, $z\in\De$ we have $[f;p/q]_\zi(z)\in\C$. Therefore, $B(f,\zi)(z)\neq0$, where $B$ is given by the Jacobi formula.

So there is a $\de''>0$ such that $\de''<1$ and $\mid B(f,\zi)(z)\mid>\de''$ for all $\zi\in L$, $z\in\De$, because $L\times\De$ is compact.

By making ``$a$'' small enough, by continuity one can get
\[
\mid B(g,\zi)(z)\mid>\frac{\de''}{2} \ \ \text{for all} \ \ \zi\in L, \ \ z\in\De.
\]
For $\el\in\{0,1,\ld,s\}$ it holds
\begin{align*}
\sup_{\zi\in L}\,\sup_{z\in\De}\mid\big[g;p/q\big]_{\zi}^{(\el)}(z)-g^{(\el)}(z)
\mid\le&\,\sup_{z\in\De}|f^{(\el)}(z)-g^{(\el)}(z)\mid \\
&+\sup_{\zi\in L}\,\sup_{z\in\de}|f^{(\el)}(z)-\big[f;p/q\big]^{(\el)}_{\zi}(z)\mid\\
&+\sup_{\zi\in L}\,\sup_{z\in\De}\mid\big[g;p/q\big]^{(\el)}_\zi(z)
-\big[f;p/q\big]^{(\el)}_\zi(z)\mid.
\end{align*}
The first term obviously get small as ``$a$'' gets small, because $L'\supset\De$. Since the second term is fixed and less than $1/s$ we must control only the last term.

But the Jacobi denominators of $\big[f;p/q\big] ^{(\el)}_\zi(z)$ and
$\big[g;p/q\big]^{(\el)}_\zi(z)$ are bounded below from
$(\de'')^{\el+1}$ and $(\de''/2)^{\el+1}$ respectively for
$\el=0,1,\ld,s$.

Thus, the last term can get as small as we want to for all $\el=0,1,\ld,s$, if $a$ is small enough.\ We are done.

Finally, we prove that for all $j,s\in\N$ the set $\bigcup\limits_{(p,q)\in F}(E(j,p,q,s)\cap T(p,q,s))$ is dense in $X^\infty(\OO)$.

Let $L'\subset\C$ be a compact set inside $\oO$ such that
$L\cup\De\subset L'$. Without loss of generality we may assume that
every connected component of $\widetilde{\C}\sm L'$ contains a point
that belongs to $\widetilde{\C}\sm\oO$.\ This can be achieved, for
example, by taking $L'=\oO\cap\overline{D(0,n)}$, for big enough
$n\in\N$.

Let also $g$ be a function inside $X^\infty(\OO)$, $N\in\N$ and $\e>0$. We can assume, without loss of generality, that $g$ is a rational function with poles off $\oO$, because of the definition of $X^\infty(\OO)$.

To prove what we want to, we have to find a function $f\in X^\infty(\OO)$ and a pair $(p,q)\in F$ such that:\vspace*{-0.2cm}
\begin{enumerate}
\item[(i)] $f\in D_{p,q}(\zi)\cap E_{p,q,\zi}(K\cup\De)$, for all $\zi\in L$.\vspace*{-0.2cm}
\item[(ii)] $\dis\sup_{\zi\in L}\,\dis\sup_{z\in K}\chi([f;p/q]_\zi(z),f_j(z))<\dfrac{1}{s}$.\vspace*{-0.2cm}
\item[(iii)] $\dis\sup_{\zi\in L}\,\dis\sup_{z\in\De}\mid f^{(\el)}(z)-\big[f;p/q\big]_\zi^{(\el)}(z)\mid<\dfrac{1}{s}$, for $\el=0,1,\ld,s$\vspace*{-0.2cm}
\item[(iv)] $\dis\sup_{z\in L'}\mid f^{(m)}(z)-g^{(m)}(z)\mid<\e$, for $m=0,1,\ld,N$.
\end{enumerate}

Let $\oo:L'\cup K\ra\C$, such that $\oo(z)=\left\{\begin{array}{cc}
                                                    f_j(z), & z\in K \\
                                                    g(z), & z\in L'
                                                  \end{array}\right.$.

Now, let $\mi$ be the sum of the principal parts of the poles of the
rational function $f_j$ that belong to $K$. Then $(\oo-\mi)$ is
holomorphic in a neighborhood of $(L'\cup K)$. Combining Runge's
with Weierstrass Theorems we conclude that there exists a rational
function $\dfrac{\tA(z)}{\tB(z)}$ with poles out of $(L'\cup K)$,
approximating $(\oo-\mi)$ uniformly on $L'\cup K$ with respect to
the euclidean metric and in the level of all derivatives of order
from zero to $N$. That implies that the function
$\dfrac{A(z)}{B(z)}=\mi(z)+\dfrac{\tA(z)}{\tB (z)}$ approximates
$f_j(z)$, uniformly on $K$ with respect to the chordal distance, and
also that $\Big(\dfrac{A(z)}{B(z)}\Big)^{(\el)}$ approximates the
function $(g(z))^{(\el)}$ uniformly on $L'$, with respect to the
euclidean metric.\ Obviously, we can assume that the greatest common
divisor of $A(z)$ and $B(z)$ is equal to one.

From our assumption on $F$, there exists a pair $(p,q)\in F$ such
that $p>\deg A,\deg B$ and $q>\deg B$. We consider the function
$\dfrac{A(z)}{B(z)}+dz^T=\dfrac{A(z)+d\cdot z^T\cdot B(z)}{B(z)}$
where $T=p-\deg B$ and $d$ is different than zero. Now, it is easy
to see that $gcd(A(z)+dz^TB(z),B(z))$ equals again to one.\ Thus,
according to Proposition \ref{prop2.2} it holds that for all
$\zi\in\C$ such that $B(\zi)\neq0$ the rational function
$\dfrac{A(z)+dz^TB(z)}{B(z)}$ belongs to $D_{p,q}(\zi)$ and also
$\Big[\dfrac{A(z)+dz^TB(z)}{B(z)};p/q\big]_\zi(z)=\dfrac{A(z)+dz^TB(z)}
{B(z)}$. In particular the above hold for all $\zi\in L'$, because
$B(\zi)\neq0$ for all $\zi\in L'$.

We distinguish the cases $B(z)\neq0$ for all $z\in\oO$ and the case
where $B$ has roots in $\oO\sm L'$. First assume that $B(z)\neq0$
for all $z\in\oO$. In this case we set
$f(z)=\dfrac{A(z)}{B(z)}+dz^T$, and by selecting $d$ with $|d|$
small enough, we are done.

In the second case, since every component of $\widetilde{\C}\sm L'$
contains a point from $\widetilde{\C}\sm\oO$, there exists a
rational function that belongs to $X^\infty(\OO)$, call it $f$, such
that every finite set of derivatives $f^{(\el)}$ are close to
$\Big(\dfrac{A(z)+dz^T}{B(z)}\Big)^{(\el)}$ uniformly on $L'$. This
is immediate from Runge's and Weierstrass Theorems and also from the
fact that $B$ has finitely many roots outside $L'$ and thus in a
positive distance from $L'$.

It is easy to see that $f$ fulfills all requirements in the same way
as $\dfrac{A(z)+dz^TB(z)}{B(z)}$ does except from the fact that
maybe $[f;p/q]_\zi(z)\neq f(z)$.

But the following is true:
\begin{align*}
\dis\sup_{\zi\in L}\,\sup_{z\in\De}\mid\big[f;p/q\big]^{(\el)}_\zi(z)-f^{(\el)}(z)\mid\le&\,\sup_{z\in\De}
\mid f^{(\el)}(z)-h^{(\el)}(z)\mid \\
&+\sup_{\zi\in L}\,\sup_{z\in\De}\mid h^{(\el)}(z)-\big[h;p/q\big]^{(\el)}_\zi(z)\mid\\
&+\sup_{\zi\in L}\,\sup_{z\in\De}\mid\big[h;p/q\big]^{(\el)}_\zi(z)-
\big[f;p/q\big]^{(\el)}_\zi(z)\mid.   \hspace*{0.8cm} \mbox{($\ast$)}
\end{align*}
with $h(z)=\dfrac{A(z)+dz^TB(z)}{B(z)}$.

Now, as $p,q$ are fixed and we can control any finite set of
derivatives of $f$, we can also control any finite set of Taylor
coefficients of $f$. Thus, we can make the first and the last term
of the right-hand side expression in ($\ast$) as small as we want to
and we are done.

This completes the proof of the Theorem. \qb
\end{Proof}

If we set $K=\emptyset$ and $L=\De=L_n$ where
$L_n=\oO\cap\overline{D(0,n)}$ for $n=1,2,\ld$ and we apply Baire's
Theorem once more, we obtain the result that generically all $f\in
X^\infty(\OO)$ can be approximated by their Pad\'{e} approximants
$[f;p_n/q_n]_\zi(z)$, $(p_n,q_n)\in F$, provided that
$F\subset\N\times\N$ contains a sequence $(\tp_n,\tq_n)\in F$,
$n=1,2,\ld$ such that $\tp_n\ra+\infty$ and $\tq_n\ra+\infty$.

If we set $L=\{\zi\}\subset\oO$, $K=K_n,\De=L_n$ where $K_n$ is an
exhausting sequence of $\C\sm\oO$ \cite{Rud} and
$L_n=\oO\cap\overline{D(0,n)}$ for $n=1,2,\ld$ and then apply
Baire's Theorem, provided that the set $F\subset\N\times\N$ contains
a sequence $((\tp_n\tq_n))_{n\in\N}$ where $\tp_n\ra+\infty$ and
$\tq_n\ra+\infty$ we obtain the following:
\begin{thm}\label{3.2}
Let $\OO\subset\C$ be an open set and $\zi\in\oO$ be fixed. Then
there exist $f\in X^\infty(\OO)$ such that, for every rational
function $h$ and every compact set $K\subset\C\sm\oO$ there exists a
sequence $(p_n,q_n)\in F$, $n=1,2,\ld$ such that $f\in
D_{p_n.q_n}(\zi)$ for all $n\in\N$, and $\dis\sup_{z\in
K}\chi([f;p_n/q_n]_\zi(z),h(z))\ra0$ and for every compact set
$L'\subset\oO$ there is a $n(L')\in\N$ such that $f\in
E_{p_n,q_n,\zi}(K\cup L')$, for all $n\ge n(L')$ and also
$\dis\sup_{z\in
L'}\mid\big[f;p_n/q_n\big]^{(\el)}_\zi(z)-f^{(\el)}(z)\ra0$, as
$n\ra+\infty$, for all $\el\in\N$.\ The set of such functions $f\in
X^\infty(\OO)$ is dense and $G_\de$ in $X^\infty(\OO)$.
\end{thm}

If we set $L=\De=L_n$ and $K=K_n$ for $n=1,2,\ld$ where $L_n=\oO\cap\overline{D(0,n)}$ and $\big(K_n\big)^\infty_{n=1}$ is an exhausting sequence of $\C\sm\oO$ by applying Baire Theorem we obtain the following:
\begin{thm}\label{thm3.3}
Let $F$ be a subset of $\N\times\N$ containing a sequence
$(\tp_n,\tq_n)\in F$, $n=1,2,\ld$ with $\tp_n\ra+\infty$,
$\tq_n\ra+\infty$.\ Let $\OO\subset\C$ be an open set.\ Then there
exists a function $f\in X^\infty(\OO)$ that satisfy the following:

For every compact set $K\subset\C\sm\oO$ and rational function $h$ there exists a sequence $(p_n,q_n)\in F$, $n=1,2,\ld$ such that the following hold:

For every compact set $L\subset\oO$ there exists a $n(L)\in\N$ such that
\[
f\in D_{p_n,q_n}(\zi), \ \ \text{for all} \ \ n\ge n(L) \ \ \text{and} \ \ \zi\in L
\]
\[
f\in E_{p_n,q_n,\zi}(K\cup L), \ \ \text{for all} \ \ n\ge n(L) \ \ \text{and} \ \ \zi\in L
\]
and
\[
\sup_{\zi\in L}\,\sup_{z\in K}\chi\big(\big[f;p_n/q_n\big]_\zi(z),h(z)\big)\ra0, \ \ \text{as} \ \ n\ra+\infty
\]
\[
\sup_{\zi\in L}\,\sup_{z\in L}\mid\big[ f;p_n/q_n\big]^{(\el)}_\zi(z)-f^{(\el)}(z)\mid\ra0, \ \ \text{as} \ \ n\ra+\infty \ \ \text{for all} \ \ \el=0,1,2,\ld.
\]
The set of such functions $f\in X^\infty(\OO)$ is dense and $G_\de$ in $X^\infty(\OO)$.
\end{thm}
\section{The case $\bbb{\{\infty\}\cup(\C\sm\oO)}$ connected}\label{sec4}
\noindent

We recall that if $K\subseteq\C$ is compact then $A(K)=\{h:K\ra\C$
continuous on $K$ and holomorphic in $K^0\}$.
\begin{thm}\label{thm4.1}
Let $F\subset\N\times\N$ be a set that contains a sequence $(\tp_n\tq_n)$, $n=1,2,\ld$ such that $\tp_n\ra+\infty$ and $\OO\subseteq\C$ an open set such that $\{\infty\}\cup(\C\sm\oO)$ is connected.

Let $L,\De\subset\C$ compact sets inside $\oO$ and $K$ be a compact
set in $\C$ such that $K^c$ is connected and $K\cap\oO=\emptyset$.
Then there exist $f\in X^\infty(\OO)$ such that:

For every function $h$ in $A(K)$ there exists a sequence $(p_n,q_n)\in F$, $n=1,2,\ld$ such that:
\begin{enumerate}
\item[(i)] $f\in D_{p,q}(\zi)\cap E_{p,q,\zi}(K\cup\De)$, for all $\zi\in L$
\item[(ii)] for all $\el=0,1,2,\ld$, $\dis\sup_{\zi\in L}\,\dis\sup_{z\in\De}\mid\big[f;p_n/q_n\big]^{(\el)}_\zi(z)-f^{(\el)}(z)\mid\ra0$
as $n\ra+\infty$
\item[(iii)] $\dis\sup_{\zi\in L}\,\dis\sup_{z\in K}\mid[f;p_n/q_n]_\zi(z)-h(z)\mid\ra0$, as $n\ra+\infty$.
\end{enumerate}

The set of such functions $f\in X^\infty(\OO)$ is dense and $G_\de$ is $C^\infty(\OO)$.
\end{thm}
\begin{Proof}
Let $\big(f_j\big)^\infty_{j=1}$ be an enumeration of all polynomial
functions is with coefficients from $\Q+i\Q$.

We name $\cu$ the set of the functions with the properties (i),
(ii), (iii) and we will prove that $\cu$ is $G_\de$-dense set in the
$X^\infty(\OO)$-topology and so $\cu\neq\emptyset$.

For $j,s\in\N^\ast$ and $(p,q)\in F$ we define:
\begin{align*}
E(j,p,q,s)=&\,\big\{f\in X^\infty(\OO)\big|f\in D_{p,q}(\zi)\cap E_{p,q,\zi}(K), \ \ \text{for all}\ \ \zi\in L \ \ \text{and} \\
&\sup_{\zi\in L}\,\sup_{z\in K}\big|[f;p/q]_\zi(z)-f_j(z)\big|<\frac{1}{s}\big\}
\end{align*}
\begin{align*}
T(p,q,s)=&\,\big\{f\in X^\infty(\OO)\big|f\in D_{p,q}(\zi)\cap E_{p,q,\zi}(\De), \ \ \text{for all} \ \ \zi\in L \ \ \text{and} \\
&\sup_{\zi\in L}\,\sup_{z\in\De}\big|\big[f;p/q\big]_\zi^{(\el)}(z)-f^{(\el)}(z)\big|<
\frac{1}{s} \ \ \text{for} \ \ \el=0,1,\ld,s\big\}.
\end{align*}
It is true that
$\cu=\bigcap\limits^\infty_{j,s=1}\bigcup\limits_{(p,q)\in
F}(E(j,p,q,s)\cap T(p,q,s))$. This can easily be verified using
Mergelyans' Theorem. Now to prove that $\cu$ is a $G_\de$-dense set
in $X^\infty(\OO)$, it is enough to show that for every
$j,s=1,2,\ld$ and $(p,q)\in F$ the sets $E(j,p,q,s)$ and $T(p,q,s)$
are open in $X^\infty(\OO)$ and that for every $j,s=1,2,\ld$, the
set $\bigcup\limits_{(p,q)\in F}(E(j,p,q,s)\cap T(p.q.s))$ is dense
in $X^\infty(\OO)$.

So, let $j,s\in\N^\ast$ and a pair $(p,q)\in F$. We first prove that
$E(j,p,q,s)$ is open in $X^\infty(\OO)$.

Indeed, let $f\in E(j,p,q,s)$, and let $g\in X^\infty(\OO)$ be such
that
\begin{eqnarray}
\sup_{z\in L}\mid f^{(m)}(z)-g^{(m)}(z)\mid<a \ \  \text{for} \ \  m=0,1,2,\ld,p+q+1  \label{eq1}
\end{eqnarray}

The number $a>0$ will be determined later on. It is enough to prove
that if $a$ is small enough then $g\in E(j,p,q,s)$.

The Hankel determinants defining $D_{p,q}(\zi)$ for $f$ depend
continuously on $\zi\in L$; thus, there exists $\de>0$ such that the
absolute values of the corresponding Hankel determinants are greater
than $\de>0$, for every $\zi\in L$, because $f\in D_{p,q}(\zi)$, for
every $\zi\in L$ and because $L$ is compact.

From (\ref{eq1}) we can control the first $p+q+1$ Taylor
coefficients of $g$ and by making $a>0$ small enough one can get the
Hankel determinants that define $D_{p,q}(\zi)$ to have absolute
value at least $\de/2>0$.

Therefore, $g$ will belong in $D_{p,q}(\zi)$ for every $\zi\in L$.
Now we consider the Pad\'{e} approximants of $f,g$ according to the
Jacobi formula (see preliminaries)
\[
[f,p/q]_\zi(z)=\frac{A(f,\zi)(z)}{B(f,\zi)(z)} \ \ \text{and} \ \ [g,p/q]_\zi(z)=\frac{A(g,\zi)(z)}{B(g,\zi)(z)}.
\]
Now $\mid A(f,\zi)(z)\mid^2+\mid B(f,\zi)(z)\mid^2$ vary
continuously with respect to $(z,\zi)\in K\times L$, because of the
Jacobi formula. So, there is a $\de'>0$ such that:
\[
\mid A(f,\zi)(z)\mid^2+\mid B(f,\zi)(z)\mid^2\ge\de', \ \ \text{for all} \ \ \zi\in L \ \ \text{and} \ \ z\in K.
\]
Now again from the Jacobi formula, if $a$ is small enough, one gets:
\[
\mid A(g,\zi)(z)\mid^2+\mid B(g,\zi)(z)\mid^2\ge\de'/2, \ \ \text{for all} \ \ \zi\in L \ \ \text{and} \ \ z\in K.
\]
This yields that $g\in E_{p,q,\zi}(K)$ for every $\zi\in L$.

Now because $\big|[f;p/q]_\zi(z)-f_j(z)\big|<\dfrac{1}{s}$, for all
$\zi\in L,z\in K$ it follows that $[f;p/q]_\zi(z)\in\C$, for all
$\zi\in L$, $z\in K$.

Thus, there exist $\de''>0$ such that $|B(f,\zi)(z)|>\de''$ for all
$\zi\in L$, $z\in K$, because $L\times K$ is compact. Because the
first $p+q+1$ Taylor coefficients of $g$ can be controlled and
because of the Jacobi formula, by making ``$a$'' small enough, one
gets
\[
\mid B(g,\zi)(z)\mid>\frac{\de''}{2}, \ \ \text{for all} \ \ \zi\in L,\ \ z\in K.
\]
To complete the proof that $E(j,p,q,s)$ is open, it is enough to
show that $\dis\sup_{\zi\in L}\,\dis\sup_{z\in K}\mid
[g;r/q]_\zi(z)-[f;p/q]_\zi(z)\mid$ can become less than
$\dfrac{1}{s}-\dis\sup_{\zi\in L}\,\dis\sup_{z\in
K}\mid[f;p/q]_\zi(z)-f_j(z)\mid=\ga>0$.

But
\[
\sup_{\zi\in L}\,\sup_{z\in K}\mid[f;p/q]_\zi(z)-[g;p/q]_\zi(z)\mid\le\frac{2}{(\de'')^2}\mid
A(f,\zi)B(g,\zi)(z)-B(f,\zi)(z)A(g,\zi)(z)\mid.
\]
This easily yields the result because the expression on the right-hand side of the inequality can become as small as we want to, for $\el=0,1,\ld,s$.

The proof that $T(p,q,s)$ is open in $X^\infty(\OO)$ is similar to the corresponding proof in Theorem \ref{thm3.1} and is omitted.

Finally, we prove that for every $j,s\in\N^\ast$ the set
$\bigcup\limits_{(p,q)\in F}(E(j,p,q,s)\cap T(p,q,s))$ is dense in
$X^\infty(\OO)$. Let $L'\subset\C$ be a compact set inside $\oO$,
such that $L\cup\De\subset L'$. We can assume, without loss of
generality, that $L'=\oO\cap\overline{D(0,n)}$ for some $n\in\N$,
big enough. Let $g$ be a function inside $X^\infty(\OO)$, $N\in\N$
and $\e>0$. We can assume by the definition of $X^\infty(\OO)$ and
from the fact that $\{\infty\}\cup(\C\sm\oO)$ is connected, that $g$
is a polynomial. We have to find a function $f$ inside
$X^\infty(\OO)$ and a pair $(p,q)\in F$ such that:\vspace*{-0.2cm}
\begin{enumerate}
\item[(i)] $f\in D_{p,q}(\zi)\cap E_{p,q,\zi}(K\cup\De)$, for all $\zi\in L$\vspace*{-0.2cm}
\item[(ii)] $\dis\sup_{\zi\in L}\,\dis\sup_{z\in \De}\big|\big[f;p/q\big]
    ^{(\el)}_\zi(z)-f^{(\el)}(z)\big|<\dfrac{1}{s}$, for all $\el=0,1,\ld,s$\vspace*{-0.2cm}
\item[(iii)] $\dis\sup_{\zi\in L}\,\dis\sup_{z\in K}\mid[f;p/q]_\zi(z)-f_j(z)\mid<\dfrac{1}{s}$ \vspace*{-0.2cm}
\item[(iv)]  $\dis\sup_{z\in L'}\mid f^{(m)}(z)-g^{(m)}(z)\mid<\e$ for $m=0,1,\ld,N$.
\end{enumerate}

Now, $\widetilde{\C}\sm
L'=\widetilde{\C}\sm(\oO\cap\overline{D(0,n)})=(\{\infty\}\cup(\C\sm\oO))\cup
(\widetilde{\C}\sm\overline{D(0,n)})$ is a connected set as a union
of two connected subsets of $\widetilde{\C}$, intersecting at least
at the point $\infty$.

From our hypothesis $K^c$ is connected too. So, as $L\cap K=\emptyset$, there exist two simply connected domains $G_1$, $G_2$ such that $G_1\cap G_2=\emptyset$, $L\subset G_1$, $K\subset G_2$. We may assume also that $G_1,G_2$ have positive distance.

We consider now the function $w:G_1\cup G_2\ra\C$ with:
$w(z)=\left\{\begin{array}{cc}
f_j(z),& z\in G_2 \\
g(z), & z\in G_1
\end{array}\right.$ . By Runge's theorem there exists a sequence of polynomials $\tp_n$ that approximate uniformly on compact sets the analytic function $w$.

Because $G_1\cup G_2$ is open, according to Weierstrass theorem the
approximation will be valid in the level of all derivatives.
Therefore, one such polynomial $\tp$ approximates $f_j$ uniformly on
$K$ with respect to the euclidean distance and $\tp ^{(\el)}$
approximate $g^{(\el)}$ with respect to the euclidean metric,
uniformly on $L'$ for all $\el=0,1,\ld,N$.

Now, there exists $(p,q)\in F$ with $p>\deg\tp$, $q\ge0$ and because
$\deg(\tp(z)+dz^p)=p$, for all $d>0$ by Proposition \ref{prop2.2} we
have $\tp(z)+dz^p\in D_{p,q}(\zi)$ and
$[\tp(z)+dz^p;p/q]_\zi(z)=\tp(z)+dz^p$ for all $\zi\in\C$.\ But
$\tp(z)+dz^p$ approximate, as $d\ra0$, the polynomial $\tp(z)$
uniformly for any finite set of derivatives and on any compact
subset of $\C$.

Therefore, if we choose $d$ sufficiently small and set
$f(z)=\tp(z)+dz^p$, we are\linebreak done. \qb
\end{Proof}

Varying $L$, $\De$ and $K$ we can obtain more complete versions of Theorem \ref{thm4.1} as we do in Section \ref{sec3} for Theorem \ref{thm3.1}.
\section{Density of rational functions}\label{sec5}
\noindent

In this section we give sufficient conditions so that $X^\infty(\OO)=A^\infty(\OO)$.
\begin{thm}\label{thm5.1}
Let $\OO$ be a bounded, connected and open set such that:\vspace*{-0.2cm}
\begin{enumerate}
\item[(a)] $(\oO)^0=\OO$.\vspace*{-0.2cm}
\item[(b)] $\C\sm\oO$ is connected. \vspace*{-0.2cm}
\item[(c)] There exists $M<+\infty$, such that for every $a,b\in\OO$ there exists a curve $\ga$ inside $\OO$ (i.e.\ $\ga:[0,1]\ra\OO$ continuous function)
such that $\ga(0)=a$, $\ga(1)=b$ and Length$(\ga)\le M$. \vspace*{-0.2cm}
\end{enumerate}

Then the polynomials are dense in $A^\infty(\OO)$ (and therefore
$X^\infty(\OO)=A^\infty(\OO)$).
\end{thm}
\begin{Proof}
Let $f\in A^\infty(\OO)$, $\e>0$ and $n\in\N_0=\{0,1,2,\ld\}$.

It suffices to find a polynomial $p$ such that:
\[
\sup_{z\in\OO}\mid f^{(\el)}(z)-p^{(\el)}(z)\mid<\e, \ \ \text{for} \ \ \el=0,1,\ld,n.
\]
Now $f^{(n)}\in C(\oO)$ and is analytic in $\OO=(\oO)^0$, because $f\in A^\infty(\OO)$.\ Also $\oO$ is a compact set as $\OO$ is bounded.

Thus, by Mergelyans' Theorem there exists a polynomial, $p_n$ such that $\dis\sup_{z\in\OO}\mid f^{(n)}(z)-p_n(z)\mid<\dfrac{\e}{(M+1)^n}$.

Now, fix $z_0\in\OO$.\ Then, for every $z\in\OO$, there exists a
curve $\ga_z$ inside $\OO$ that starts at $z_0$ and ends at $z$ and
has length at most $M$.

Also, for $0\le k\le n-1$, we define the polynomial $p_k(z)$ by:
\[
p_k(z)=f^{(k)}(z_0)+\int_{[z_0,z]}p_{k+1}(\zi)d\zi
\]
and we set $p=p_0$.\ Then it is obvious that
$p^{(k)}=p^{(k)}_0=p_k$, for all $k$, $0\le k\le n$.

Now for $k=n$ we have:
\begin{align*}
\sup_{z\in\OO}\mid f^{(k)}(z)-p^{(k)}(z)\mid&=\sup_{z\in\OO}\mid f^{(n)}
(z)-p_n(z)\mid<\frac{\e}{(M+1)^n}=\frac{e}{(M+1)^k}. \\
&\text{Therefore,}\;\sup_{z\in\OO}\mid f^{(k)}(z)-p^{(k)}(z)\mid<\frac{\e}{(M+1)^k}.  \hspace*{2.3cm} \text{($\ast$)}
\end{align*}
Assume that the above relationship holds for a fixed $k$, $1\le k\le
n$. We will prove it for $k-1$: It is:
\begin{align*}
\sup_{z\in\OO}\mid f^{(k-1)}(z)-p^{(k-1)}(z)\mid&=\sup_{z\in\OO}\bigg|\int_{\ga_z}(f^{(k)}(\zi)-p^{(k)}(\zi))d\zi\bigg| \\
&\le\sup_{z\in\OO}\int_{\ga_z}\mid f^{(k)}(\zi)-p^{(k)}(\zi)|\,|d\zi|\\
&\le\sup_{z\in\OO}\int_{\ga_z}\frac{\e}{(M+1)^k}\mid d\zi\mid\le\frac{M\e}{(M+1)^k}<\frac{\e}{(M+1)^{k-1}}.
\end{align*}
which is exactly what we wanted.

That means that ($\ast$) is true for all $k=0,1,\ld,n$ and our proof is complete. \qb
\end{Proof}
\begin{thm}\label{thm5.2}
Let $\OO$ be a connected, open set such that \vspace*{-0.2cm}
\begin{enumerate}
\item[(a)] $(\oO)^0=\OO$.\vspace*{-0.2cm}
\item[(b)] $\{\infty\}\cup(\C\sm\oO)$ is a connected set \vspace*{-0.2cm}
\item[(c)] There exists $n_0\in\N$ such that for every $N\ge n_0$, there exists $M_N>0$ such that for all $a,b\in\overline{\OO\cap D(0,N)}^0$ there exists a continuous function $\ga:[0,1]\ra\overline{\OO\cap D(0,N)}^0$ with $\ga(0)=a$, $\ga(1)=b$ and $Length(\ga)\le M_N$.\vspace*{-0.2cm}
\end{enumerate}

Then the polynomials are dense in $A^\infty(\OO)$, and therefore
$X^\infty(\OO)=A^\infty(\OO)$.
\end{thm}

For the proof we need two lemmas.
\begin{lem}\label{lem5.3}
Let $\OO$ be an open set in $\C$, such that
$\{\infty\}\cup(\C\sm\oO)$ is connected. Then, for every $N\in\N$,
$\{\infty\}\cup(\C\sm\overline{(D(0,N)\cap\OO)})$ is connected.
\end{lem}
\noindent
{\bf Proof of Lemma \ref{lem5.3}}.
Let $V$ be a connected component of $\{\infty\}\cup(\C\sm\overline{(D(0,N)\cap\OO)}$). It is enough to prove that $\infty\in V$.

Now, $\{\infty\}\cup(\C\sm\overline{(D(0,N)\cap\OO)})$ is a non
empty open set inside $\C\cup\{\infty\}$ and therefore, $V$ is a non
empty open set in $\C\cup\{\infty\}$.

That implies that there exist $x\in V$ with $|x|<N$, or there exists
$x\in V$ with $|x|>N$.

In the former case, $x\in D(0,N)$ and as $x\notin\overline{(\OO\cap D(0,N)})$, one can see that $x\notin\oO$.\ Indeed, if not there would be $(x)_{n\in\N}\subset\OO$ with $x_n\ra x$.\ But $D(0,N)$ is an open set. Thus, eventually, we have $x_n\in D(0,N)$. This implies that for $n$ big enough we have $x_n\in (\OO\cap D(0,N))$.\ It follows that $x\in\overline{(\OO\cap D(0,N))}$, which is a contradiction.

Thus, $x\in(\C\sm\oO)\cup\{\infty\}\subseteq\{\infty\}\cup(\C\sm\overline{D(0,N)\cap\OO})$.

Because $\{\infty\}\cup(\C\sm\oO)$ is connected, it follows that $\infty\in V$ as we wanted.

In the latter case, we have
\[
x\in\{\infty\}\cup(\C\sm\overline{D(0,N)})\subseteq\{\infty\}\cup
(\C\sm\overline{(\OO\cap D(0,N)}).
\]
Because $\{\infty\}\cup(\C\sm\overline{D(0,N)}$ is connected, it follows again this $\infty\in V$.

This completes the proof of the Lemma \ref{lem5.3}.  \qb
\begin{lem}\label{lem5.4}
Let $\OO$ be an open set, $\OO\subseteq\C$.\ Then for every $N\in\N$ we have
\[
\overline{\overline{(\OO\cap D(0,N))}^0}=\overline{(\OO\cap D(0,N))}.
\]
\end{lem}
\noindent
{\bf Proof of the Lemma \ref{lem5.4}}.\! Since $\overline{\OO\cap
D(0,N)}\!\supseteq\!\OO\cap D(0,N)$ it follows $\overline{(\OO\cap
D(0,N))}^0\!\supseteq\!(\OO\cap D(0,N))^0$.

This implies $\overline{(\OO\cap D(0,N))}^0\supseteq(\OO\cap D(0,N))$, as the latter set is open.

Therefore, $\overline{\overline{(\OO\cap D(0,N))}^0}\supseteq\overline{(\OO\cap D(0,N))}$.

For the other inclusion let $x\in\overline{\overline{(\OO\cap D(0,N))}^0}$.\  Then there exist $x_n\in\overline{(\OO\cap D(0,N))}^0$ such that $x_n\ra x$.\ Now, for every $n\in\N$, there exist $\e_n\in\Big(0,\dfrac{1}{n}\Big)$ such that $B(x_n,\e_n)\subseteq\overline{(\OO\cap D(0,N))}$. As $x_n\in\overline{(\OO\cap D(0,N))}$, for every $n\in\N$, there exists $y_n\in B(x_n,\e_n)$ with $y_n\in\OO\cap D(0,N)$. But $s_n\ra0$ and $x_n\ra x$ which gives $y_n\ra x$. Thus, $x\in\overline{(\OO\cap D(0,N))}$ and the proof of the Lemma \ref{lem5.4} is completed. \qb\vspace*{0.2cm} \\
%
%
{\bf Proof of the Theorem \ref{thm5.2}}. Let $f\in A^\infty(\OO)$, $\e>0$, $n\in\N_0$ and $N\in\N$.\ It is enough to find a polynomial $p$, such that $\dis\sup_{z\in\overline{(\OO\cap D(0,N))}}\mid f^{(\el)}(z)-p^{(\el)}(z)\mid<\e$, for $\el=0,1,\ld,n$.

Without loss of generality we can assume, $N\ge n_0$. Now let $V=(\overline{(\OO\cap D(0,N))})^0$.

$V$ is an open, connected and bounded set.\ ($V$ is connected
because $N\ge n_0$ and because of condition (c) of our hypothesis).

Also, from Lemma \ref{lem5.4}, we get
$\oV=\overline{\overline{(\OO\cap D(0,N))}^0}=\overline{\OO\cap
D(0,N))}$ and therefore, $\oV^0=\overline{(\OO\cap D(0,N))}^0=V$.

From Lemmas \ref{lem5.3} and \ref{lem5.4} it follows that $\{\infty\}\cup(\C\sm\oV)=\{\infty\}\cup(\C\sm(\overline{\OO\cap D(0,N)})$ is connected.

Since $V$ is bounded, it follows that $C\sm\oV$ is connected.

Thus, $V$ satisfies all conditions of Theorem \ref{thm5.1} and
therefore the set of all polynomials is dense in $A^\infty(V)$. But
$f\in A^\infty(\OO)$ and $\OO\supset V$. Indeed we have
$\OO\supset\OO\cap D(0,N)$.\ This gives
$\oO\supset\overline{(\OO\cap D(0,N))}$.\ Therefore,
$\OO=\oO^0\supset( \overline{\OO\cap D(0,N)})^0=V$.\ Thus,
$\OO\supset V$.

But $\OO\supset V$ implies $A^\infty(\OO)\subset A^\infty(V)$.

Thus, $f\in A^\infty(V)$.\ Therefore, there exists a polynomial $p$
such that
\[
\sup_{z\in V}\mid f^{(\el)}(z)-p^{(\el)}(z)\mid<\frac{\e}{2}, \ \ \text{for} \ \ \el=0,1,\ld,n.
\]
This implies
\[
\sup_{z\in \oV}\mid f^{(\el)}(z)-p^{(\el)}(z)\mid<\frac{\e}{2}<\e, \ \ \text{for} \ \ \el=0,1,\ld,n.
\]
By Lemma \ref{lem5.4} we have that
\[
\oV=\overline{(\overline{(\OO\cap D(0,N)))}^0)}=\overline{\OO\cap D(0,N)}.
\]
The proof of Theorem \ref{thm5.2} is complete. \qb
\begin{thm}\label{thm5.5}
Let $\OO$ be a bounded, connected, open set such that:\vspace*{-0.2cm}
\begin{enumerate}
\item[(a)] $(\oO)^0=\OO$.\vspace*{-0.2cm}
\item[(b)] $\{\infty\}\cup(\C\sm\oO)$ has exactly $k$ connected components, $k\in\N$.\vspace*{-0.2cm}
\item[(c)] There exists $M>0$ such that for all $a,b\in\OO$, there exists a continuous function $\ga:[0,1]\ra\OO$ with $\ga(0)=a$, $\ga(1)=b$ and $Length(\ga)\le M$.\vspace*{-0.2cm}
\end{enumerate}

Now pick from each connected component of $\{\infty\}\cup(\C\sm\oO)$
a point $a_i$, $i=0,1,2,\ld,\linebreak k-1$ and set
$S=\{a_0,\ld,a_{k-1}\}$, where $a_0$ belongs to the unbounded
component. Then the set of all rational functions with poles from
$S$ is dense in $A^\infty(\OO)$ and therefore
$X^\infty(\OO)=A^\infty(\OO)$.
\end{thm}

For the proof we need the following lemma.
\begin{lem}\label{lem5.6}
Let $\OO$ be an open set, $n\in\N$, $n\ge1$, and let $f$ be
holomorphic in $\OO$. Then, for $i\in\N$, $0\le i\le n-1$, and $\ga$
any closed curve in $\OO$ of bounded variation we have
$\int\limits_\ga z^if^{(n)}(z)dz=0$.
\end{lem}
{\bf Proof of Lemma \ref{lem5.6}}. We use induction on $n$. For
$n=1$, we have to prove $\int\limits_\ga f'(z)dz=0$, which is
obvious, because the curve $\ga$ is closed.\ Suppose the statement
holds for $n=k$, we will prove it for $k+1$. Let
$i\in\{1,\ld,(k+1)-1=k\}$, we have
\setcounter{equation}{0}
\begin{eqnarray}
\int_\ga z^if^{(k+1)}(z)dz=\int_\ga z^i(f^{(k)}(z)'dz=z^if^{(k)}(z)\bigg|^{\ga(1)}_{\ga(0)}-i\int_\ga z^{i-1}f^{(k)}(z)dz . \label{eq1}
\end{eqnarray}
Since $\ga$ is a closed curve it follows that
\begin{eqnarray}
z^if^{(k)}(z)\bigg|_{\ga(0)}^{\ga(1)}=0.  \label{eq2}
\end{eqnarray}
Since $i-1<k$ from the induction hypothesis we have
\begin{eqnarray}
\int_\ga z^{i-1}f^{(k)}(z)dz=0.  \label{eq3}
\end{eqnarray}
Relations (\ref{eq1}), (\ref{eq2}) and (\ref{eq3}) imply
\[
\int_\ga z^if^{(k+1)}(z)dz=0,
\]
as we wanted.

Finally, for $i=0$ we have $\int\limits_\ga
z^0f^{(k+1)}(z)dz=\int\limits_\ga(f^{(k)}(z))'dz=0$, because $\ga$
is closed.

The proof of the lemma is complete. \qb \vspace*{0.2cm} \\
{\bf Proof of Theorem \ref{thm5.5}}. We first prove it in the case $a_0=\infty$.

Let $f\in A^\infty(\OO)$, $\e>0$ and $n\in\N=\{0,1,2,\ld\}$.\ We
need to find a rational function $r$ with poles in $S$ such that:
\[
\sup_{z\in\OO}\mid f^{(\el)}(w)-r^{(\el)}(w)|<\e, \ \ \text{for} \ \ \el=0,1,\ld,n.
\]
Fix curves $\ga_i$, $i=1,2,\ld,k-1$ in $\OO$ which are closed, have finite length and also $\text{Ind}(\ga_i,a_j)=\left\{\begin{array}{cc}
                                          1, & i=j \\
                                          0, & i\neq j
                                        \end{array}\right.$. For the construction of such curves see \cite{Ahlfors}.
Since $f\in A^\infty(\OO)$, it follows that $f^{(n)}\in C(\oO)$ and that $f^{(n)}$ is analytic in $\OO=(\oO)^0$.

Thus, from Mergelyans' Theorem (\cite{Rud}, ch. 20, ex 1) there
exists a rational function called $\tr_n(z)$, with poles only in
$S$, such that $\dis\sup_{z\in\OO}\mid f^{(n)}(z)-\tr_n(z)\mid<a$
where $a>0$ is sufficiently small.

In particular it suffices that
\[
0<a<\min\bigg(\frac{\e}{2(M+1)^n},\frac{\e\cdot\pi\cdot \ti^n}{n(k-1)\cdot\Big(D+\ssum^{k-1}_{i=1}\mid a_i\mid\Big)^n}\cdot
\frac{1}{(M+1)^n\Big(\ssum^{k-1}_{i-1}length(\ga_i)\Big)}\bigg)
\]
where $D\ge\max(1,\dis\max_{z\in\oO}\mid z\mid)$ and
$0<\ti\le\min(1,\dist(S,\OO))$.

Now, by analyzing $\tr_n(z)$ into simple fractions, there exists a rational function $r_n(z)$ with poles only in $S$ such that
\[
\tr_n(z)=r_n(z)+\sum^{k-1}_{i=1}\sum^n_{j=1}\frac{b_{ij}}{(z-a_i)^j} \ \ \text{with} \ \ b_{ij}\in\C \ \ i=1,2,\ld,k-1 \ \ \text{and} \ \ j=1,2,\ld,n.
\]
and
\[
Res((z-a_i)^{j-1}r_n(z),a_i)=0, \ \ \text{for all} \ \ i=1,2,\ld,k-1, \ \ j=1,2,\ld,n.
\]
Fix $(i,j)\in\{1,2,\ld,k-1\}\times\{1,2,\ld,n\}$.\ Then using Lemma \ref{lem5.6}, it follows that
\begin{align*}
\mid b_{ij}\mid&=\bigg|\frac{1}{2\pi i}\int_{\ga_i}(z-a_i)^{j-1}\tr_n(z)dz\bigg|=\bigg|\frac{1}{2\pi i}\int_{\ga_i}
(z-a_i)^{j-1}(\tr_n(z)-f^{(n)}(z))dz\bigg| \\
&\le\frac{1}{2\pi}(D+|a_i|)^nL(\ga_i)\cdot a.
\end{align*}
Now, choosing the positive number $a$ sufficiently small, we get
\setcounter{equation}{0}
\begin{eqnarray}
\mid b_{ij}\mid\le\frac{\e\cdot \ti^n}{2n(K-1)\cdot(M+1)^n}, \ \ \text{for all} \ \ i=1,2,\ld,K-1 \ \ \text{and} \ \ j=1,2,\ld,n.  \label{eq1}
\end{eqnarray}
Since $\dis\sup_{z\in\OO}\mid f^{(n)}(z)-\tr_n(z)|<a$, it follows that
\[
\sup_{z\in\OO}\bigg(\mid f^{(n)}(z)-r_n(z)\mid-\sum^{k-1}_{i=1}\sum^n_{j=1}\frac{\mid b_{ij}\mid}{\mid z-a_i\mid^j}\bigg)<a.
\]
This implies,
\[
\sup_{z\in\OO}\mid f^{(n)}(z)-r_n(z)\mid<a+\sum^{k-1}_{i-1}\sum^n_{j=1}\frac{\mid b_{ij}\mid}{\mid z-a_i\mid^j}\le\frac{\e}{2(M+1)^n}+
\sum^{k-1}_{i=1}\sum^n_{j=1}\frac{\mid b_{ij}\mid}{\ti^n}
\]
Combining this with relation (\ref{eq1}) we obtain,
\begin{align*}
\sup_{z\in\OO}\mid f^{(n)}(z)-r_n(z)\mid&\le\frac{\e}{2(M+1)^n}+\sum^{k-1}_{i=1}\sum^n_{j=1}
\frac{\e\cdot \ti^{n-j}}{2n(k-1)\cdot(M+1)^n}\\
&\le\frac{\e}{2(M+1)^n}+\frac{\e}{2(M+1)^n}=\frac{\e}{(M+1)^n}.
\end{align*}
The function $r_n$ has a Laurent expansion around each $a_i\in S\sm\{\infty\}$, where the coefficients of $(z-a)^\el$ for $\el=-n,-n+1,\ld,-1$ are equal to zero.

This implies that for each $s,1\le s\le n$ the integral
$\underset{\leftarrow\ \  s \ \
\ra}{\int\limits\int\limits\cdots\int\limits}r_n(z)(dz)^s$ defines a
regular holomorphic function in $\OO$, which is not multivalued.

We proceed by induction on $\la\in\{n,n-1,\ld,0\}$. For $\la\in\N$,
$0\le\la\le n-1$, we define:
\[
r_\la(z)=f^{(\la)}(z_0)+\int_{[z_0,z]}r_{\la+1}(z)dz,
\]
where $r_{\la+1}$ is known by the induction hypothesis.

Thus, we define the rational functions $r_n,r_{n+1},\ld,r_1,r_0$.\ We set $r=r_0$.

It is obvious that $r_\la(z)=r_0^{(\la)}(z)=r^{(\la)}(z)$, for
$\la=0,1,\ld,n$. The proof of the case $a_0=\infty$ can be completed
as the last part of the proof of Theorem \ref{thm5.1}.

Next we consider the general case where $a_0$ is not necessarily equal to $\infty$.

Let $f\in A^\infty(\OO)$, $\e>0$ and a natural $n\in\N=\{0,1,2,\ld\}$. We seek a rational function $r$ with poles only in $S=\{a_0,\ld,a_{k-1}\}$, such that $\dis\sup_{z\in\OO}\mid f^{(\el)}(z)-r^{(\el)}(z)\mid<\e$ for $\el=0,1,\ld,n$.

From the previous case, there exists a rational function $r_1$ with
poles only in $\tS=\{\infty,a_1,a_2,\ld,a_{k-1}\}$ such that
\begin{eqnarray}
\sup_{z\in\OO}\mid f^{(\el)}(z)-r^{(\el)}_1(z)\mid<\frac{\e}{2}, \ \ \text{for} \ \ \el=0,1,\ld,n.  \label{eq2}
\end{eqnarray}
But it is known that there exists a rational function $r$ with poles in $S=\{a_0,a_1,\ld,a_{n-1}\}$ such that
\begin{eqnarray}
\sup_{z\in\OO}\mid r^{(\el)}(z)-r^{(\el)}_1(z)\mid<\frac{\e}{2} \ \ \text{for} \ \ \el=0,1,\ld,n.  \label{eq3}
\end{eqnarray}
See \cite{D-M-T}, Lemma 2.2.

From relations (\ref{eq2}) and (\ref{eq3}) we derive that
\[
\sup_{z\in\OO}\mid r^{(\el)}(z)-f^{(\el)}(z)\mid<\e \ \ \text{for} \ \ \el=0,1,\ld,n
\]
and $r$ has its poles in $S$.

The proof of Theorem \ref{thm5.5} is complete now.
\qb\vspace*{0.2cm}

The following theorem is the more general one.
\begin{thm}\label{thm5.7}
Let $\OO$ be a connected, open set such that: \vspace*{-0.2cm}
\begin{enumerate}
\item[(a)] $(\oO)^0=\OO$.\vspace*{-0.2cm}
\item[(b)] $\{\infty\}\cup(\C\sm\oO)$ has exactly $k$ connected components, $k\in\N$.\vspace*{-0.2cm}
\item[(c)] There exists $n_0\in\N$ such that for every $n\ge n_0$,
there exists $M_n>0$ such that for all
$a,b\in(\overline{(\OO\cap D(0,n))})^0$, there exists a
continuous function $\ga:[0,1]\ra(\overline{\OO\cap D(0,n)})^0$
with $\ga(0)=a$, $\ga(1)=b$ and $Length(\ga)\le M_n$.

Now pick from each of the $k$ connected components of
$\{\infty\}\cup(\C\sm\oO)$ a point, $a_i$ $(i=0,1,\ld,k-1)$ and
set $S=\{a_0,a_1,\ld,a_{k-1}\}$.\vspace*{-0.2cm}
\end{enumerate}

Then, the set of rational functions with poles only in $S$ is dense
in $A^\infty(\OO)$, and therefore $X^\infty(\OO)=A^\infty(\OO)$.
\end{thm}
\begin{Proof}
Let $r>0$ be such that $D(0,r)$ contains all the components of
$\{\infty\}\cup(\C\sm\oO)$ not containing $\infty$.\ This is
possible, since $k\in\N$.

Let $f\in A^\infty(\OO)$, $\e>0$, $n\in\N=\{0,1,2,\ld\}$ and
$N\in\N$, $N\neq0$. It is enough to find a rational function $r$
with poles only in $S$ such that: $\dis\sup_{z\in\overline{(\OO\cap
D(0,N))}^0}\mid f^{(\el)}(z)-r^{(\el)}(z)\mid<\e$, because
$\overline{(\overline{\OO\cap D(0,N))}^0}=\overline{\OO\cap D(0,N)}$
(Lemma \ref{lem5.4}).

Without loss of generality we may assume that $N\ge n_0+r$. We claim that there exists $M>0$, such that for every $a,b\in\overline{\OO\cap D(0,N)}^0$, there exists a curve in $\overline{\OO\cap D(0,N)}^0$  that joins $a$ and $b$ and has length at most $M$ and also that the set $\{\infty\}\cup(\C\sm\overline{(\OO\cap D(0,N))}$ has at most $k$ connected components, each of them containing at least one point from $S$.

The former is immediate according to our hypothesis by setting $M=M_N$.

For the latter, let $V$ be a connected component of
$\{\infty\}\cup(\C\sm(\overline{\OO\cap D(0,N)})$. Because $V$ is
open and non empty there exists $x\in V$ with $|x|>N$ or there
exists $x\in V$ with $|x|<N$.

In the first case we have that
$x\in\{\infty\}\cup(\C\sm\overline{(D(0,N)})\subseteq\{\infty\}\cup(\C\sm
\overline{(\OO\cap D(0,N))})$. Because
$\{\infty\}\cup(\C\sm\overline{(D(0,N)})$ is connected, it follows
that $\infty\in V$. Thus, the unbounded component of
$\{\infty\}\cup(\C\sm\oO)$ is contained in $V$, which implies that
$V\cap S\neq\emptyset$.

In the latter case, $x\in(\{\infty\}\cup\C\sm\overline{(\OO\cap
D(0,N))})\cap D(0,N)$. Therefore, $x\in (\C\sm\overline{(\OO\cap
D(0,N))})\cap D(0,N)$.

It follows that $x\notin\oO$. Indeed, if not, there exists a sequence $(x_n)_{n\in\N}\subset\OO$ with $x_n\ra x$. But $D(0,N)$ is open. Thus, ther exists $n_0\in\N$, such that $x_n\in(\OO\cap D(0,N))$ for every $n\ge n_0$.

It follows that $x\in\overline{(\OO\cap D(0,N))}$, contradicting the assumption that $x$ belongs to $\C\sm\overline{(\OO\cap D(0,N))}$.

Therefore $x\in\C\sm\oO\subseteq\{\infty\}\cup(\C\sm\oO)\subseteq\{\infty\}\cup
(\C\sm\overline{(\OO\cap D(0,N))})$.

Let $V_1$ be the connected component inside
$\{\infty\}\cup(\C\sm\oO)$ containing $x$. It follows that
$V_1\subset V$ and thus, $V\cap S\neq\emptyset$.

Therefore, we have proved that any connected component of
$\{\infty\}\cup(\C\sm\overline{(\OO\!\cap\! D(0,N)})$ intersets non
trivially $S$. Since $S$ contains exactly $k$ points and the
components of $\{\infty\}\cup(\C\sm\overline{(\OO\cap D(0,N)})$ are
mutually disjoint we conclude that the number of the components is
at most $k$. (It can also be proved that, if $N$ is big enough, the
number of components is exactly $k$, but this is not needed at the
sequel).

Now, set $T=\overline{(\OO\cap D(0,N))}^0$. We can easily check that $T$ satisfies all assumptions of Theorem \ref{thm5.5}. It is true that $f\in A^\infty(\OO)$.\ But $T\subset\oO^0=\OO$. Thus, $A^\infty(T)\supset A^\infty(\OO)$. This implies that $f\in A^\infty(T)$. Theorem \ref{thm5.5}, combined with the fact that the set $S$ contains at least one point from each component of $\{\infty\}\cup(\C\sm\tiT)$ implies that, there exists a rational function $r$ with poles only in $S$ such that,
\[
\sup_{z\in\bar{T}}|f^{(\el)}(z)-r^{(\el)}(z)|<\e, \ \ \text{for} \ \  \el=0,1,2,\ld,n.
\]
It follows that,
\[
\sup_{z\in\overline{(\OO\cap D(0,N))}}|f^{(\el)}(z)-r^{(\el)}(z)|<\e, \ \ \text{for} \ \  \el=0,1,2,\ld,n,
\]
because $\bar{T}=\overline{\OO\cap D(0,N)}$, according to Lemma
\ref{lem5.4}. The proof is complete. \qb
\end{Proof}
\begin{rem}\label{rem5.8}
Lemma \ref{lem5.6} can be generalized to a necessary and sufficient condition for an analytic function to have an antiderivative in $\OO$ of order $n$, $n\in\N$.

More specifically it holds the following:

Let $n\in\N$, $\OO$ an open subset of $\C$ and $f$ an analytic function in $\OO$. The following are equivalent \vspace*{-0.2cm}
\begin{enumerate}
\item[(a)] There exists a function $F$, which is analytic in $\OO$, such that $F^{(n)}(z)=f(z)$, for all $z\in\OO$.\vspace*{-0.2cm}
\item[(b)] For any closed curve $\ga$ in $\OO$ of bounded variation, and for every $i=0,1,\ld,n-1$, it is true that $\int\limits_\ga z^if(z)dz=0$.\vspace*{-0.2cm}
\item[(c)] For any closed curve $\ga$ in $\OO$ of bounded variation, and for every polynomial $P$ with $\deg P\le n-1$, it is true that $\int\limits_\ga P(z)f(z)dz=0$.\vspace*{-0.2cm}
\end{enumerate}
\end{rem}
\begin{rem}\label{rem5.9}
The previous theorems remain valid for any open set $\OO\subseteq\C$ with a finite number of components, where each component satisfies the preresquities of the according theorem under the extra condition that the closures of any two components are disjoint.
\end{rem}

The proof of this is immediate by applying our theorems in each
component. This gives a finite number of rational functions, one for
each component.

Applying Runge and Weierstrass theorems we find one rational function with poles off $\oO$ approximating simultaneously the above rational functions.

We do not have the answer in the case where the closure of the
components are not disjoint but we know that the answer is positive
in the particular case of two open discs $D_1,D_2$ such that
$\oD_1\cap\oD_2$ is a singleton. Indeed, let $f\in A^\infty(D_1\cup
D_2)$, $\e>0$ and $n\in\N$.\ We can assume that the disks touch at
zero and also that their radius is at most 1. The open set $D_1\cup
D_2$ obviously satisfies the \vspace*{-0.2cm}
\begin{enumerate}
\item[(i)] $\overline{(D_1\cup D_2)}^0=D_1\cup D_2$ \vspace*{-0.2cm}
\item[(ii)] $\C\sm\overline{(D_1\cup D_2)}$ is connected and also \vspace*{-0.2cm}
\item[(iii)] for any two points in $D_1\cup D_2$, there exists a
polygonal line joining them of length at most four, that may be
 touches the boundary at most at zero, and otherwise is
 contained in $D_1\cup D_2$. \vspace*{-0.2cm}
\end{enumerate}

Now, for $f\in A^\infty(D_1\cup D_2)$ and $\ga$ a polygonal simple
curve in $\overline{D_1\cup D_2}$, that touches the boundary at most
at zero, it is immediate from an argument of continuity that
\[
\int_\ga f^{(\la)}(z)dz=f^{(\la)}(\ga(1))-f^{(\la)}(\ga(0)), \ \ \text{for any} \ \ \la\in\N.
\]
The rest of the proof is similar to the proof of Theorem \ref{thm5.1}.
\begin{rem}\label{rem5.10}
It can be shown that if $\OO$ is a Jordan domain with rectifiable
boundary, it fulfills the preresquities of Theorem \ref{thm5.1}.

More specifically it holds that there exists a positive constant
$M>0$ such that any two points in $\OO$ can be joined by a curve
inside $\OO$ of length at most $M$.

Moreover the above holds in the case of a domain bounded by $k$
disjoint Jordan curves with rectifiable boundaries.

Indeed in the case of a Jordan domain $\OO$ with rectifiable
boundary (as in \cite{M-N}) every point in $\OO$ is joined with the
boundary with a segment with length at most $\text{diam}(\OO)$. Next
two points on the boundary of $\OO$ can by joined by subarc of the
boundary with length at most the length of the boundary. Thus
$M=2diam(\OO)+length(\partial\OO)$.

However, the curve is not contained in $\OO$. According to a Theorem
of Caratheodory \cite{Koosis} every conformal map $\f$ from the open
unit disc $D$ onto $\OO$ extends to a homeomorphism from $\oD$ to
$\oO$.\ Further since the boundary of $\OO$ is rectifiable, it
follows that $\f'\in H^1$ \cite{Koosis}, \cite{Duren}. Thus we can
use the image $\Ga$ by $\f$ of a circumference $C(0,r)$, $0<r<1$,
where $r$ is very close to 1 and we can replace the subarc of
$\partial\OO$ by an arc of $\Ga$; its length is less than or equal
to $\|\f'\|_1$ which is equal to the length of $\partial\OO$.

When we have $k$ disjoint Jordan curves with rectifiable boundaries,
first we join the outer boundary with another boundary using a
segment of minimum length (which is minimum for all boundaries).
This segment is disjoint from all other boundaries.\ Let $E_1$, be
the compact set containing the two previous boundaries and the
segment.\ We joint $E_1$ with some other boundary using a segment of
minimum length.\ We continue in this way and after a finite number
of steps we obtain a (connected) curve $E$ containing all boundaries
and whose all other points belong to $\OO$.\ The length of $E$ is
finite. If we consider to points $z_1,z_2\in\OO$ we join each one of
them with some boundary using two segments.\ Then we joint $z_1,z_2$
by these two segments and a piece of $E$.\ The length does not
exceed $2\text{diam}(\OO)+length(E)$.\ For the arcs contained in the
outer boundary wee can use the conformal mapping $\f:D\ra\OO'$,
where $\OO'$ is the union of $\OO$ with all bounded component of
$\C\sm\oO$, which is simply connected.

Thus, without increasing the length we can have a curve joining
$z_1,z_2$ in $\oO$ not meeting the outer boundary. For another
boundary $\ga_j$ let $b$ be a point interior to $\ga_j$.\ Thus,
$\dist(b,\oO)=r>0$. Using the inversion $w=\dfrac{1}{z-b}$, the
complement of the interior of $\ga_j$ (with $\infty$ included) is
transformed to a bounded simply connected domain $\OO''$ containing
0.\ Using again a conformal mapping $g:D\ra\OO''$ (and $g'$ is again
in $H^1$) we can replace the subarc of $\ga_j$ contained in our
curve by another arc inside $\OO$. Its length may be increased but
we can have it as close to the initial length as we wish.\ Thus, the
assumptions of Theorem \ref{thm5.5} are satisfied with
$M=2diam(\OO)+length(E)+\de$, for any $\de>0$.\ In particular we can
have
\[
M=2diam\OO+length(E)+1.
\]
\end{rem}
\begin{rem}\label{rem5.11}
We can have examples of Jordan domains $\OO$ without rectifiable
boundary, but satisfying the assumptions of Theorem \ref{thm5.1}.
For instance, if we consider any starlike Jordan domain $\OO$ then
any two points may be joined in $\OO$ be a curve consisting of two
segments and therefore its length does not exceed $2diam(\OO)$.
Certainly we can arrange that the boundary of $\OO$ has infinite
length.

Further another example is the following.

Let $\f:[0,1]\ra\R$ be continuous, $c<\min\{\f(t):t\in[0,1]\}$ and
$\OO=\{(x,y):0<x<1$, $c<y<\f(x)\}$. If $\f$ is not of bounded
variation then the length of $\partial\OO$ is infinite but the
assumptions of Theorem \ref{thm5.2} are satisfied.
\end{rem}
\begin{rem}\label{rem5.12}
An alternative proof of Theorem \ref{thm5.5} is by using the
statement of Theorem \ref{thm5.1} combined by the following Laurent
decomposition (\cite{Co-N-P}).

Let $\OO$ be a domain of finite connectivity. Let
$V_0,V_1,\ld,V_\el$ be the components of $(\C\cup\{\infty\})\sm\OO$,
where $\infty\in V_0$. Let $f\in A^\infty(\OO)$, then
$f=f_0+f_1+\cdots+f_\el$ where $f_j\in A^\infty[(V^c_j)^0]$ for
$j=0,1,\ld,\el$ and $\dis\lim_{z\ra\infty}f_j(z)=0$ for
$j=1,\ld,\el$.
\end{rem}
\begin{rem}\label{rem5.13}
If $\OO$ is a domain satisfying the assumptions of Theorem
\ref{thm5.1} or Theorem \ref{thm5.2} or Theorem \ref{thm5.5} or
Theorem \ref{thm5.7} or the stronger assumptions discussed in this
section, then $X^\infty(\OO)=A^\infty(\OO)$.\ Therefore in these
cases the results of the Section \ref{sec3} and \ref{sec4} become
generic in $A^\infty(\OO)$.
\end{rem}
\section{Smoothness of the integration operator}\label{sec6}
\noindent

It is known that if $D=\{z\in\C\mid\,|z|<1\}$ and $f\in H^1$ then
$F(z)=\int\limits^z_0f(\zi)d\zi$ has an absolute convergent Taylor
series in $\oD$ (Hardy Inequality).\ \cite{Duren}, \cite{Hoffman}
Thus, in $D$ the antiderivative of a bounded analytic function is
also a bounded analytic function.

Moreover, in $D$ the antiderivative of a function in $A(D)$ is also contained in $A(D)$.

Now, if $\OO$ is a Jordan domain, a theorem of Caratheodory states
that every Riemann conformal mapping $\phi:D\ra\OO$ extends to an
homeomorphism $\phi:\oD\ra\oO$.\ \cite{Koosis} Moreover, the
boundary of $\OO$ is rectifiable if and only if $\phi'\in H^1$
\cite{Duren},\cite{Koosis}.

Combining the statements above, we see that if $\OO$ is a Jordan
domain with rectifiable boundary then the antiderivative of any
bounded analytic function defined on $\OO$, is also a bounded
analytic function on $\OO$.\ Furthermore, the antiderivative can be
extended continuously to $\oO$.\ More specifically the
antiderivative of every function in $A(\OO)$ remains again in
$A(\OO)$.

We will now examine the case where the Jordan domain $\OO$ does not have rectifiable boundary.
\begin{prop}\label{prop6.1}
Let $\OO$ be a Jordan domain such that there exist a constant $M<+\infty$ with the property that any two points inside $\OO$ can be joined with a curve inside $\OO$ of length at most $M$.\ Let $f$ be a bounded analytic function on $\OO$; then the antiderivative of $f$ is also a bounded analytic function on $\OO$.
\end{prop}
\begin{Proof}
Fix $z_0\in\OO$, and for every $z\in\OO$, choose a curve $\ga_z$ in $\OO$ that joins $z_0$ and $z$ and has length at most $M$. Then the antiderivative $F(z)$ is equal to $\int\limits_{\ga_z}f(\zi)d\zi$ and the result easily follows. \qb
\end{Proof}

It is easy to find examples of Jordan domains $\OO$ with non
rectifiable boundary that satisfy the preresquities of Proposition
\ref{prop6.1}, as we discussed in Remark \ref{rem5.11}.

For example, consider a starlike domain with no rectifiable boundary or the case of a domain
\[
\OO=\{(x,y)\mid0<x<1,\;c<y<\ti(x)\}, \ \ \text{where} \ \ \ti:[0,1]\ra\R
\]
is a continuous function with no bounded variation and $c<\dis\min_{x\in[0,1]}\ti(x)$.\ We call the last domain ``Domain of type $\ast$''.

Furthermore we have
\begin{prop}\label{prop6.2}
Let $\OO$ be a starlike domain or a domain of type $\ast$. Let $f\in A(\OO)$; then the antiderivative of $f$ belongs also to $A(\OO)$.
\end{prop}
\begin{Proof}
We give the proof only in the case of a Jordan domain $\OO$ which is starlike; the proof in the case of a domain of type $\ast$ is similar and is omitted.

Assume that $\OO$ is a bounded domain which is starlike with respect
to a point $z_0\in\OO$, say $z_0=0$. If $f\in A(\OO)$, it follows
that $f$ is uniformly continuous. Thus, if $\e_1>0$ is given, there
exists $\de>0$, $\de<\e_1$, so that $|f(P)-f(Q)|<\e_1$ for all
$P,Q\in\OO$ with $|P-Q|<\de$. One antiderivative of $F$ is given by
\[
F(z)=\int_{[0,z]}f(\zi)d\zi=\int^1_0f(tz)\cdot zdt \ \ \text{for} \ \ z\in\OO.
\]
It suffices to show that $F$ is uniformly continuous on $\OO$ and therefore $F\in A(\OO)$. If $z,w\in\OO$ are such that $|z-w|<\de$, it follows that $\mid tz-tw\mid<\de$ for all $t\in[0,1]$.\ Therefore, $|f(tz)-f(tw)|<\e_1$ and $\mid z-w\mid<\de$.\ It follows that $|f(tz)z-f(tw)w|\le|f(tz)|\cdot|z-w|+|f(tz)-f(tw)|\,|w|\le\|f\|_\infty\cdot\de+\e_1diam(\OO)\le\e_1[\|f\|_\infty+diam(\OO)]<\e$
provided that $\e_1$ has been chosen so that $0<\e_1<\dfrac{\e}{\|f\|_\infty+diam(\OO)}$.\ This completes the proof. \qb
\end{Proof}

After these statements, it is natural to ask whether there exists a
Jordan domain $\OO$ and a function $f\in A(\OO)$ such that the
antiderivative of $f$ is not in $A(\OO)$.

We provide such a counter-example finding a Jordan domain $\OO$ and a function $f\in A(\OO)$ such that the integral of $f$ is not even bounded inside $\OO$.
\begin{prop}\label{prop6.3}
There exist a Jordan domain $\OO$ and a function $f\in A(\OO)$ such
that the antiderivative of $f$ is not bounded inside $\OO$.
\end{prop}
\begin{Proof}
Consider the function $g:\C\ra\C$, defined by
$g(z)=(z-1)\exp\Big(\dfrac{z+1}{z-1}\Big)$ for $z\neq1$ and
$g(1)=0$.\ According to Proposition \ref{prop2.3} there exists a
Jordan domain, in the upper half plane that contains an arc of the
unit circle, having 1 as one of its endpoints such that $g$ is one
to one there. Call this Jordan domain $V$, and set $\OO=g(V)$.

Define $f:\OO=g(V)\ra\C$ by $f(w)=\dfrac{1}{\log(1-g^{-1}(w))}\cdot\dfrac{1}{\exp\Big(\dfrac{g^{-1}(w)+1}{g^{-1}(w)-1}\Big)}$. It is easy to see that $f\in A(\OO)$.

Now, consider points $z_0,z$ in the unit circle and in $V$. If the
antiderivative of $f$ was bounded, then
$\Big|\int\limits^{g(z)}_{g(z_0)}f(\zi)d\zi\Big|\le M$, for every
$z$ in the unit circle and in $V$, for some constant $M<+\infty$.

The above gives $\Big|\int\limits^z_{z_0}f(g(\zi))g'(\zi)d\zi\Big|\le M$ for every $z$ in the unit circle close enough to 1, from the upper half plane.

Thus, setting $z_0=e^{t_0}$ and $z=e^{i\tit}$, with $t$ close to
$0^+$, we have
\[
\bigg|\int^{\tit}_{t_0}f(g(e^{it}))g'(e^{it})\cdot e^{it}dt\bigg|\le M, \ \ \text{for every} \ \ 0<\tit<t_0.
\]
A computation gives
\[
f(g(e^{it}))g'(e^{it})e^{it}=f(g(e^{it}))\exp\bigg(\frac{e^{it}+1}{e^{it}-1}\bigg)
\frac{e^{it}-3}{e^{it}-1}=\frac{e^{it}-3}{e^{it}-1}\cdot\frac{1}{\log(1-e^{it})}\cdot e^{it}.
\]
Thus it must holds that
\setcounter{equation}{0}
\begin{eqnarray}
\bigg|\int^{\tit}_{t_0}e^{it}\frac{(e^{it}-3)}{(e^{it}-1)}\cdot
\frac{1}{\log(1-e^{it})}dt\bigg|\le M  \label{eq1}
\end{eqnarray}
for every $\tit$, with $0<\tit<t_0$. But if we set
$h(t)=e^{it}\Big(\dfrac{e^{it}-3}{e^{it}-1}\Big)\dfrac{1}{\log(1-e^{it})}$,
then we have that $h$ satisfies the preresquities of the Lemma
\ref{lem2.4}.

Indeed,the M\"{o}bius function $z\ra\dfrac{z-3}{z-1}$, sends the
unit circle to the line: $\{z\in\C|\Re(z)=2\}$ and also satisfies
the fact that
$\dis\lim_{t\ra0^+}\arg\Big(\dfrac{e^{it}-3}{e^{it}-1}\Big)=\dfrac{\pi}{2}$.
It is also obvious that
\[
\lim_{t\ra0^+}\arg(e^{it})=0 \ \  \text{and} \ \ \lim_{t\ra0^+}\arg\bigg(\frac{1}{\log(1-e^{it})}\bigg)=-\lim_{t\ra0^+}\arg(\log(1-e^{it}))=\pi.
\]
Thus the $\dis\lim_{t\ra0^+}\arg(h(t))$ exists. According to Lemma \ref{lem2.4}, one can easily check that the integral $\int\limits^{\tit}_{0^+}|h(t)|dt$ has the same nature with $\int\limits^{\tit}_{0^+}\dfrac{1}{t|\ln t|}dt=+\infty$.

This means that $\Big|\int\limits^{\tit}_{t_0}h(t)dt\Big|$ cannot be
bounded for all $\tit:0<\tit<t_0$,yielding\linebreak the desired
contradiction with relation 1.\ The proof of Proposition
\ref{prop6.3} is now\linebreak complete. \qb
\end{Proof}

This conversation leads us to Volterra operators on
$D=\{z\in\C\mid\,|z|<1\}$. Let $g$ be an analytic function on $D$.\
Then the operator $T_g$ maps an analytic function $f$ on $D$ to the
antiderivative of $fg'$ vanishing at 0.

Open problems in this area are to characterize the functions $g$
such that for all $f\in H^\infty(D)$, it holds that $T_g(f)\in
H^\infty(D)$, and also to characterize the functions $g$, such that
for all $f\in A(D)$, it holds $T_g(f)\in A(D)$ see \cite{Hawi}. It
is obvious that if $g'\in H^1$ then both are satisfied.\ If $\OO$ is
a starlike domain or a domain of type $\ast$ without rectifiable
boundary, then Proposition \ref{prop6.2} yields for the Riemann
mapping $g:D\ra\OO$, that $T_g(H^\infty(D))\subset H^\infty(D)$
despite the fact that $g'\notin H^1$.

Moreover, it also holds $T_g(A(D)\subset A(D)$.

In a more general way, if for a Jordan domain it holds that there
exists $M<+\infty$, such that every two points in $\OO$ can be
joined by a curve inside $\OO$ with length at most $M$, then the
Riemann mapping $g:D\ra\OO$ satisfies $T_g(H^\infty(D))\subset
H^\infty(D)$, according to Proposition \ref{prop6.1}. This happens,
more specifically, even if the boundary of $\OO$ is not rectifiable.
\section{Some generic results}\label{sec7}
\noindent

In the case of the Jordan domain $\OO$ constructed in Proposition \ref{prop6.3} the set of functions $f\in A(\OO)$, such that their antiderivative $F$ is not bounded is not void and in fact it is $G_\de$ and dense in $A(\OO)$.\ This follows from the following proposition.
\begin{prop}\label{prop7.1}
Let $\OO$ be a Jordan domain in $\C$.\ We consider the sets.\vspace*{0.2cm} \\
$X_1(\OO)=\{f\in A(\OO)$: The antiderivative of $f$ does not belong to $H^\infty(\OO)\}$.\vspace*{0.2cm} \\
$X_2(\OO)=\{f\in A(\OO)$: The antiderivative of $f$ does not belong to $A(\OO)\}$.\vspace*{0.2cm} \\
$X_3(\OO)=\{f\in H^\infty(\OO)$: The antiderivative of $f$ does not belong to $H^\infty(\OO)\}$ and\vspace*{0.2cm} \\
$X_4(\OO)=\{f\in H^\infty(\OO)$: The antiderivative of $f$ does not
belong to $A(\OO)\}$.

Then we have

i) If $X_1(\oo)\neq\emptyset$, then $X_1(\OO)$ is dense and $G_\de$ in $A(\OO)$ and $X_2(\OO)$ is residual in $A(\OO)$ and

ii) If $X_3(\OO)\neq\emptyset$, then $X_3(\OO)$ is dense and $G_\de$ in $H^\infty(\OO)$ and $X_4(\OO)$ is residual in $H^\infty(\OO)$.
\end{prop}
Proposition \ref{prop7.1} follows easily from the following.
\begin{prop}\label{prop7.2}
For $g\in H(D)$ we consider the Volterra operator $T_g:H(D)\ra
H(D)$, where $T_g(f)$ is the antiderivative of $fg'$ vanishing at 0
for any $f\in H(D)$.\ We consider the following sets:
\begin{align*}
Y_1(g)&=\{f\in A(D):\,T_g(f)\notin H^\infty(D)\} \\
Y_2(g)&=\{f\in A(D):\,T_g(f)\notin A(D)\} \\
Y_3(g)&=\{f\in H^\infty(D):\,T_g(f)\notin H^\infty(D)\} \ \ \text{and} \\
Y_4(g)&=\{f\in H^\infty(D):\,T_g(f)\notin A(D)\}.
\end{align*}
Then we have

i) If $Y_1(g)\neq\emptyset$, then $Y_1(g)$ is dense and $G_\de$ in $A(D)$ and $Y_2(g)$ is residual in $A(D)$.

ii) If $Y_3(g)\neq\emptyset$, then $Y_3(g)$ is dense and $G_\de$ in $H^\infty(D)$ and $Y_4(g)$ is residual in $H^\infty(D)$.
\end{prop}
\begin{Proof}
We have the following description of
$Y_1(g):\,Y_1(g)=\bigcap\limits_{M\in \N}E_M(g)$, where
$E_M(g)=\{f\in A(D):\|T_g(f)\|_\infty>M\}$.

First we will show that $E_M(g)$ is open in $A(\OO)$; this will
imply that $Y_1(g)$ is $G_\de$ in $A(\OO)$. Let $f\in E_M(g)$; then
there exists $z_0;|z_0|<1$, such that
$\Big|\int\limits^{z_0}_0f(\zi)g'(\zi)d\zi|>M$.\ Let $\e_1>0$ to be
defined later.\ If $\tf\in A(\OO)$ satisfies $\|f-\tf\|_\infty<\e_1$
then
$|T_g(\tf)(z_0)-T_g(f)(z_0)|\le\e_1\dis\sup_{|\zi|\le|z_0|}|g'(\zi)|\cdot|
z_0|\le\e_1\dis\sup_{|\zi|\le|z_0|}|g'(\zi)|$. Thus,
$|T_g(\tf)(z_0)|\ge|T_g(f)(z_0)|-\e_1\dis\sup_{|\zi|\le|z_0|}|g'(\zi)|$.

We choose $\e_1>0$, such that,
$|T_g(f)(z_0)|-\e_1\dis\sup_{|\zi|\le|z_0|}|g'(\zi)|>M$. It
follows\linebreak $\|T_g(\tf)\|_\infty>M$ and $\tf\in E_M(g)$.
Therefore $E_M(g)$ is open in $A(\OO)$.

Next we show that $E_M(g)$ is dense in $A(\OO)$.\ If we do so then Baire's Category Theorem will complete the proof. Let $w\in A(D)$ and $\e>0$; we are looking for a function $f\in A(\OO)$ such that $\|w-f\|_\infty<\e$ and $\|T_g(f)\|_\infty>M$.

Since $Y_1(g)\neq\emptyset$, there exists a function $\el\in A(\OO)$ such that $\|T_g(\el\|_\infty=+\infty$.\ It suffices to set $f=w+\e_1\el$ where $\e_1>0$ is sufficiently small. Then $\|T_g(f)\|=+\infty>M$.

It follows that $E_M(g)$ is $G_\de$ dense in $A(\OO)$. Since $Y_1(g)\subseteq Y_2(g)$ the proof of i) is complete.

The proof of ii) is similar and is omitted. \qb
\end{Proof}

Next we have the following.
\begin{prop}\label{prop7.3}
Let $X=H(D)$ endowed with the topology of uniform convergence on each compact subset of $D$. Or $X=A(D)$ endowed with the supremum norm.

Then the sets $\{g\in X:\,Y_1(g)\neq\emptyset\}\equiv L_1(X)$ and $\{g\in X:\,Y_3(g)\neq\emptyset\}\equiv L_2(X)$ are dense in $X$.
\end{prop}
\begin{Proof}
We consider the Jordan domain $\OO$ of Proposition \ref{prop6.3} and
let $g_0:D\ra\OO$ be a Riemann map of $D$ onto $\OO$. Then $g_0\in
A(D)\subset X$ and $Y_1(g_0)\neq\emptyset$.\ Let $f_0\in
A(D):\,T_{g_0}(f_0)\notin H^\infty(D)$. We also have
$Y_3(g_0)\neq\emptyset$.

We will show that $L_1(X)$ is dense in $X$. Let $\oo\in X$.\ If
$Y_1(\oo)\neq\emptyset$ then $\oo\in
L_1(X)\subset\overline{L_1(X)}$. Suppose $Y_1(\oo)=\emptyset$.\ Then
$T_\oo(f_0)\in H^\infty(D)$. It follows that $T_{\oo+\e
g_0}(f_0)=T_\oo(f_0)+\e T_{g_0}(f_0)\notin H^\infty(D)$ for all
$\e>0$. Since $\dis\lim_{\e\ra0}\oo+\e g_0=\oo$ and $\oo+\e g_0\in
L_1(X)$, it follows $\oo\in\overline{L_1(X)}$.\ Thus, $L_1(X)$ is
dense in $X$.\ The proof that $L_2(X)$ is dense in $X$ is similar.
\qb
\end{Proof}

In Proposition \ref{prop7.3} we wonder if $L_1(X)$ and $L_2(X)$ are also $G_\de$ in $X$. We also wonder if we can find a complete metric topology in the set of all Jordan domains (contained in a closed disc), so that generally for all such Jordan domains $\OO$ the result of Proposition \ref{prop6.3} holds. \vspace*{0.2cm} \\
{\bf Acknowledgement.} We would like to thank Professors P. M.
Gauthier and M. Papadimitrakis for helpful communications.
 \vspace*{1cm}
 University of Athens \\
 Department of Mathematics  \\
 157 84 Panepistemiopolis \\
 Athens \\
 GREECE\medskip\\
 e-mail addresses: vnestor@math.uoa.gr \\
\hspace*{2.9cm} ilias\_\_91@hotmail.com

\end{document}